\documentclass[12pt, a4paper]{amsart}
\usepackage{amscd, amssymb, pb-diagram}
\usepackage{palatino,euler,color}
\usepackage{amsthm}
\usepackage[all]{xy}
\allowdisplaybreaks

\def\Aut{\operatorname{Aut}}

\def\newspan{\operatorname{span}}
\def\clsp{\overline{\operatorname{span}}}

\def\ker{\operatorname{ker}}

\def\id{\operatorname{id}}

\def\Ind{\operatorname{Ind}}
\def\ind{\operatorname{ind}}

\def\max{\operatorname{max}}

\def\id{\operatorname{id}}

\def\tr{\operatorname{tr}}

\def\C{\mathbb{C}}

\def\R{\mathbb{R}}
\def\N{\mathbb{N}}
\def\Z{\mathbb{Z}}
\def\T{\mathbb{T}}

\def\Q{\mathbb{Q}}

\def\TT{\mathcal{T}}
\def\LL{\mathcal{L}}
\def\OO{\mathcal{O}}
\def\KK{\mathcal{K}}

\def\HH{\mathcal{H}}
\def\AA{\mathcal{A}}

\def\QQ{\mathcal{Q}}
\def\VV{\mathcal{V}}

\def\NT{\mathcal{N}\mathcal{T}}

\newcommand{\nx}{\mathbb N^{\times}}

\newcommand{\nxnx}{{\mathbb N \rtimes \mathbb N^\times}}
\newcommand{\qxqx}{{\mathbb Q \rtimes \mathbb Q^*_+}}
\newcommand{\qn}{\mathcal Q_\mathbb N}

\newcommand{\primes}{\mathcal P}

\def\add{\textup{add}}
\def\mult{\textup{mult}}

\newtheorem{thm}{Theorem}[section]
\newtheorem{cor}[thm]{Corollary}

\newtheorem{lemma}[thm]{Lemma}
\newtheorem{prop}[thm]{Proposition}

\theoremstyle{definition}

\theoremstyle{remark}
\newtheorem{remark}[thm]{Remark}

\newtheorem{example}[thm]{Example}

\numberwithin{equation}{section}

\begin{document}

\title[$C^*$-algebras associated to dilation matrices]{\boldmath{Phase transition on Exel crossed products\\ associated to dilation matrices}}
\author{Marcelo Laca}

\author{Iain Raeburn}

\author{Jacqui Ramagge}

\address{Marcelo Laca, Department of Mathematics and Statistics\\
University of Victoria\\
Victoria, BC V8W 3P4\\
Canada}
\email{laca@math.uvic.ca}

\address{Iain Raeburn, Department of Mathematics and Statistics\\
University of Otago\\
PO Box 56\\Dunedin 9054\\
New Zealand}
\email{iraeburn@maths.otago.ac.nz}

\address{Jacqui Ramagge, School of Mathematics and Applied Statistics\\
University of Wollongong\\
NSW 2522\\
Australia}
\email{ramagge@uow.edu.au}

\date{25 January 2011}
\thanks{This research was supported by the Natural Sciences and Engineering Research Council of Canada, by the University of Otago, and by the Australian Research Council.}

\begin{abstract}
An integer matrix $A\in M_d(\Z)$ induces a covering $\sigma_A$ of $\T^d$ and an endomorphism $\alpha_A:f\mapsto f\circ \sigma_A$ of $C(\T^d)$ for which there is a natural transfer operator $L$. In this paper, we compute the KMS states on the Exel crossed product $C(\T^d)\rtimes_{\alpha_A,L}\N$ and its Toeplitz extension. We find that $C(\T^d)\rtimes_{\alpha_A,L}\N$ has a unique KMS state, which has inverse temperature $\beta=\log|\det A|$. Its Toeplitz extension, on the other hand, exhibits a phase transition at $\beta=\log|\det A|$, and for larger $\beta$ the simplex of KMS$_\beta$ states is isomorphic to the simplex of probability measures on $\T^d$.
\end{abstract}

\maketitle

\section{Introduction}\label{sec:Exel}

Actions of the real line $\R$ on $C^*$-algebras are used to describe the time evolution in physical models, and also arise in a wide variety of mathematical contexts. The KMS states for the action were originally intended to be mathematical realisations of the equilibrium states in statistical mechanics \cite{bra-rob}. More recently, mathematicians have found actions of $\R$ on algebras of number-theoretic origin that exhibit phase transitions of the sort one might expect in a statistical-mechanical model \cite{bos-con, laca98, lr}. Here we describe a similar phenomenon for the gauge action on an Exel crossed-product $C^*$-algebra associated to an integer dilation matrix $A$.

An illuminating example for the analysis of KMS states is the action $\sigma$ lifted from the gauge action of $\T$ on the Toeplitz-Cuntz algebra $\TT\OO_n$ \cite{eva}. The system $(\TT\OO_n,\sigma)$ has a single KMS state for each inverse temperature $\beta\geq \log n$, but only the one at $\beta=\log n$ factors through the purely infinite simple quotient $\OO_n$ (see, for example, \cite[Example~2.8]{ln}) . Our situation is similar: the Exel crossed product is purely infinite simple and has a unique KMS state, which has inverse temperature $\beta=\log |\det A|$, whereas its Toeplitz analogue has KMS$_\beta$ states for all $\beta\geq \log|\det A|$. Here, though, the simplex of KMS$_\beta$ states is large for $\beta> \log|\det A|$, and we have a phase transition at $\beta= \log|\det A|$.

Before stating our results more precisely, we set up some notation. We consider  a matrix $A\in M_d(\Z)$ with nonzero determinant, and write $\sigma_A$ for the associated self-covering of $\T^d=\R^d/\Z^d$. Then $\alpha_A:f\mapsto f\circ \sigma_A$ is an endomorphism of $C(\T^d)$, 
\begin{equation}\label{defL}
L(f)(z)=\frac{1}{|\det A|}\sum_{\sigma_A(w)=z} f(w)
\end{equation}
defines a transfer operator $L$ for $\alpha_A$, and the triple $(C(\T^d),\alpha_A,L)$ is one of the Exel systems studied in \cite{ehr}. We write $M_L$ for the associated right-Hilbert bimodule over $C(\T^d)$, which has underlying space $C(\T^d)$, actions defined by $f\cdot m\cdot g=fm\alpha_A(g)$, and inner product defined by $\langle m,n\rangle=L(m^*n)$. We write $\phi$ for the homomorphism of $C(\T^d)$ into $\LL(M_L)$ which implements the left action. 

If $\Sigma$ is a set of coset representatives  for $\Z^d/A^t\Z^d$, then the characters $\{\gamma_m:z\mapsto z^{m}\;:\;m\in\Sigma\}$, viewed as continuous functions on $\T^d$ and hence as ele\-ments of $M_L$, form an orthonormal basis for $M_L$ (this observation is due to Packer and Rieffel \cite{pr}, and a proof consistent with our notation is given in \cite[Lemma~2.6]{ehr}). The reconstruction formula for this basis implies that $\phi(f)$ is the finite-rank operator $\sum_m\Theta_{f\cdot \gamma_m,\gamma_m}$ for every $f\in C(\T^d)$. Then since $\alpha_A$ is unital, the results of \cite{br} imply that 
Exel's Toeplitz algebra $\TT(C(\T^d),\alpha_A,L)$ is the Toeplitz algebra $\TT(M_L)$, and that the Exel crossed product $C(\T^d)\rtimes_{\alpha_A,L}\N$ is isomorphic to the Cuntz-Pimsner algebra $\OO(M_L)$.

The Toeplitz algebra $\TT(M)$ of a Hilbert bimodule $M$ over $C$ is generated by a universal representation $(i_{M},i_{C})$ of $M$, and carries a gauge action of $\T$ characterised by $\gamma_z(i_M(m))=zi_M(m)$ and $\gamma_z(i_C(c))=i_C(c)$; this action descends to the Cuntz-Pimsner algebra $(\OO(M),j_M,j_C)$. The gauge actions inflate to actions $\sigma$ of $\R$ which are characterised by
\begin{align}\label{defdynamics}
&\sigma_t\circ i_{C}=i_C\text{ and }\sigma_t(i_M(m))=e^{it}i_M(m),\text{ and }\\ &\sigma_t\circ j_{C}=j_C\text{ and }\sigma_t(j_M(m))=e^{it}j_M(m).\notag
\end{align} 
Our goal is the following description of the KMS states of $(\TT(M_L),\sigma)$ and $(\OO(M_L),\sigma)$.

\begin{thm}\label{summary}
Suppose that $A\in M_d(\Z)$ has nonzero determinant, $(C(\T^d),\alpha_A,L)$ is the associated Exel system, and $\sigma$ denotes the action of $\R$ on $\TT(C(\T^d),\alpha_A,L)$ satisfying \eqref{defdynamics}. 

\smallskip
\textnormal{(a)} There are no KMS$_\beta$ states on $(\TT(C(\T^d),\alpha_A,L),\sigma)$ unless $\beta\geq\log|\det A|$.

\smallskip
\textnormal{(b)} For each $\beta\in(\log|\det A|,\infty]$, the simplex of KMS$_\beta$ states on $(\TT(C(\T^d),\alpha_A,L),\sigma)$ is affinely homeomorphic to the simplex $P(\T^d)$ of probability measures on $\T^d$.

\smallskip
\textnormal{(c)} If $A$ is a dilation matrix, then $(\TT(C(\T^d),\alpha_A,L),\sigma)$ has a unique KMS$_{\log|\det A|}$ state, and this state factors through the quotient map
\[
Q:(\TT(C(\T^d),\alpha_A,L),\sigma)\to (C(\T^d)\rtimes_{\alpha,L}\N,\sigma).
\]

\textnormal{(d)} Every ground state of $(\TT(C(\T^d),\alpha_A,L),\sigma)$ is a KMS$_\infty$ state.
\end{thm}

After a short review of notation and coventions, we begin in \S\ref{sec-pres} by giving presentations of $\TT(C(\T^d),\alpha_A,L)$ and $C(\T^d)\rtimes_{\alpha,L}\N$ in terms of a unitary repesentation $u$ of $\Z^d$ and an isometry $v$ which, loosely, implements the action $\alpha_A$. Then in \S\ref{sec-charKMS}, we characterise the KMS states of $(\TT(C(\T^d),\alpha_A,L),\sigma)$ in terms of their behaviour with respect to the presentation in \S\ref{sec-pres}. 

We then set about proving Theorem~\ref{summary} in stages, and we give more precise formulations of our results as we go. For example, we prove part (c) in \S\ref{sec-Exelcp}, and we prove a little more than we stated above: we only need to assume that $A$ is a dilation matrix to get uniqueness of the KMS$_{\log|\det A|}$ state. In \S\ref{sec-existence}, we prove existence of lots of KMS states (see Proposition~\ref{existforbeta>1}); a novelty in our construction is the use of induced representations to build Hilbert spaces where we can construct KMS states from vector states. Then in \S\ref{sec-par}, we prove that we have found all the KMS$_\beta$ states for $\beta>\log|\det A|$. In \S\ref{sec-ground}, we prove part (d).

Theorem~\ref{summary} and our strategy for proving it were motivated by our previous work in \cite{lr}, or more precisely, by what it says about the KMS states of the additive boundary quotient $(\TT_{\add}(\nxnx),\sigma)$ of the Toeplitz algebra $\TT(\nxnx)$ (see \cite[\S4]{BaHLR}). The connection with Exel crossed products is made in \cite[\S5]{BaHLR}, where it is shown that there is an Exel system $(C(\T), \alpha,L,\nx)$ of the kind studied in \cite{l} whose Nica-Toeplitz crossed product $\NT(C(\T),\alpha,L,\nx)$ is $\TT_{\add}(\nxnx)$ and whose Exel crossed product $C(\T)\rtimes_{\alpha,L}\nx$ is the Crisp-Laca boundary quotient of $\TT(\nxnx)$ (or in other words, Cuntz's $\QQ_{\N}$ \cite{c2}). So our present analysis differs from that in \cite{lr} in that we have raised the dimension of the torus to $d$, but have replaced $\nx\cong\N^\infty$ by $\N$.  The case $d=1$, where $A$ has the form $(N)$, is in some sense an intersection of our results with those in \cite{lr}, and in \S\ref{BSsemigp} we carry out an analysis of the KMS states on $\TT(\N\rtimes_N\N)$ parallel to that in \cite[\S4]{BaHLR}.

We close in \S\ref{sec-gpalgs} with a discussion of the case where $A$ is invertible over $\Z$. The endomorphism $\alpha_A$ is then  an automorphism, and $C(\T^d)\rtimes_{\alpha_A,L}\N$ is the usual crossed product $C(\T^d)\rtimes_{\alpha_A}\Z$, which can also be viewed as a group algebra $C^*(\Z^d\rtimes_{A^t}\Z)$. We know from \S\ref{sec-Exelcp} that there can only be KMS$_\beta$ states when $\beta=\log|\det A|=0$, so we are left to determine the invariant traces, which we do in Proposition~\ref{invtr}. When $\Z^d\rtimes_{A^t}\Z$ is the integer Heisenberg group, for example, we can find lots of invariant traces. 

\section{Notation and coventions}

\subsection{Integer matrices.} Throughout this paper $A$ is a matrix in $M_d(\Z)$ whose determinant $\det A$ is nonzero. If the eigenvalues $\lambda$ of $A\in M_d(\Z)$ all satisfy $|\lambda|>1$, then we call $A$ a \emph{dilation matrix}. This was a standing assumption in \cite{ehr}, but here we do not in general assume that $A$ is a dilation matrix. 
We use multiindex notation, so that $e^{2\pi ix}=(e^{2\pi ix_1},\cdots,e^{2\pi i x_d})$ for $x\in \R^d$, and the covering map $\sigma_A:\T^d\to \T^d$ is characterised by $\sigma_A(e^{2\pi ix})=e^{2\pi iAx}$ for $x\in \R^d$. Since the transpose $A^t$ appears more often than $A$, we write $B:=A^t$; we have tried to avoid using the letters $A$ and $B$ for anything else. We choose a set $\Sigma$ of coset representatives for $\Z^d/B\Z^d$, and assume for convenience that $0\in\Sigma$. We sometimes write $N$ for $|\det A|=|\det B|$.

\subsection{Hilbert bimodules} A bimodule $M$ over a $C^*$-algebra $C$ is a right-Hilbert bimodule if it is a right Hilbert $C$-module, and if the left action of $C$ is implemented by a homomorphism $\phi$ of $C$ into the $C^*$-algebra $\LL(M)$ of adjointable operators. (Such bimodules are also called ``correspondences'' over $C$, or just ``Hilbert bimodules'' for short.) Our $C^*$-algebras will always have identities, and our bimodules are always essential in the sense that $\phi:C\to \LL(M)$ is unital. 

A representation\footnote{These are often called Toeplitz representations, but we now believe this to have been an unfortunate choice of name (see \cite[Remark~5.3]{BaHLR}).} $(\psi,\pi)$ of a Hilbert bimodule $M$ in a $C^*$-algebra $D$ consists of a linear map $\psi:M\to D$ and a unital representation $\pi:C\to D$ such that 
\[
\psi(c_1\cdot m\cdot c_2)=\pi(c_1)\psi(m)\pi(c_2)\ \text{ and }\ \pi(\langle m,n\rangle)=\psi(m)^*\psi(n).
\]
Every Hilbert bimodule $M$ has a Toeplitz algebra $\TT(M)$, which is generated by a universal representation $(i_M,i_C)$.

Every representation $(\psi,\pi)$ of $M$ in $D$ induces a representation $(\psi,\pi)^{(1)}:\KK(M)\to D$ such that $(\psi,\pi)^{(1)}(\Theta_{m,n})=\psi(m)\psi(n)^*$ (see \cite[page~202]{p} or \cite[Proposition~1.6]{fr}). The representation $(\psi,\pi)$ is Cuntz-Pimsner covariant if
\[
(\psi,\pi)^{(1)}(\phi(a))=\pi(a) \text{ whenever } \phi(a)\in\KK(M).
\]
The Cuntz-Pimsner algebra $\OO(M)$ is the quotient of $\TT(M)$ that is universal for Cuntz-Pimsner covariant representations. We write $Q:\TT(M)\to \OO(M)$ for the quotient map, and $(j_M,j_C):=(Q\circ i_M,Q\circ i_C)$ for the universal Cuntz-Pimser covariant representation in $\OO(M)$. (Though there are several different definitions of Cuntz-Pimsner covariance out there, they all coincide for the bimodules in this paper.)

\subsection{Exel crossed products} An Exel system consists of an endomorphism $\alpha$ of a $C^*$-algebra $C$, and a transfer operator $L$ for $\alpha$, which is a bounded positive linear map $L:C\to C$ such that $L(\alpha(c)d)=cL(d)$. The examples of interest here are the systems $(C(\T^d),\alpha_A,L)$ discussed in the introduction, where $\alpha_A$ is the endomorphism $f\mapsto f\circ\sigma_A$ associated to an integer matrix $A$, and $L$ is defined by averaging over inverse images of points, as in \eqref{defL}. Notice that both $\alpha$ and $L$ are unital.

Every Exel system $(C,\alpha,L)$ gives rise to a Hilbert bimodule over $C$ as follows. We first make a copy $C_L$ of $C$ into a bimodule over $C$ by setting $c\cdot m=cm$ and $m\cdot c=m\alpha(c)$ for $m\in C_L$ and $c\in C$. The formula $\langle m,n\rangle:=L(m^*n)$  carries a $C$-valued pre-inner product on $C_L$, and completing $C_L$ gives a right Hilbert $C$-module $M_L$. Because $L$ is bounded, the left action of $C$ extends to an action of $C$ by adjointable operators on the completion $M_L$. (The details are in \cite[\S3]{br}.)
In general the completion process involves modding out by vectors of length zero, so that the quotient carries a $C$-valued inner product. However, for the systems $(C(\T^d),\alpha_A,L)$, the module $C(\T^d)$ has no vectors of length zero and is already complete (see \cite[Lemma~3.3]{LR2}). So we dispense with the quotient maps $q:C_L\to M_L$ which were used in \cite{br} to distinguish between elements of the algebra and elements of the bimodule.

For an Exel system $(C,\alpha,L)$, we define the Toeplitz algebra $\TT(C,\alpha,L)$ to be $\TT(M_L)$, and the Exel crossed product $C(\T^d)\rtimes_{\alpha_A,L}\N$ to be the Cuntz-Pimsner algebra $\OO(M_L)$. This is not quite Exel's original definition \cite{e1}, but for the systems $(C(\T^d),\alpha_A,L)$ of interest to us it is equivalent. (The precise relationship between Exel's crossed product and Cuntz-Pimsner algebras is worked out in \cite[\S3]{br}.) So  $\TT(C(\T^d),\alpha_A,L)$ and $C(\T^d)\rtimes_{\alpha_A,L}\N$ are generated by universal representations $(i_{M_L},i_{C(\T^d)})$ and $(j_{M_L},j_{C(\T^d)})$.

\subsection{KMS states} Suppose that $\sigma$ is an action of $\R$ by automorphisms of a $C^*$-algebra $C$. An element $c$ of $C$ is analytic if $t\mapsto \sigma_t(c)$ is the restriction of an entire function. A state $\phi$ of $C$ is a KMS state at inverse temperature $\beta\in (0,\infty)$ if there is a set $S$ of analytic elements such that $\newspan S$ is dense in $C$ and
\[
\phi(dc)=\phi(c\sigma_{it}(d))\ \text{ for $c,d\in S$.}
\]
In \cite[\S7]{lr}, we were careful to explain why this definition is equivalent to that used in the standard sources \cite{bra-rob} and \cite{ped}. We also adopt two more recent conventions which are possibly nonstandard. First, we regard the KMS$_0$ states to be the \emph{$\sigma$-invariant} traces; this agrees with the convention in \cite{ped} rather than the one in \cite{bra-rob}. Second, we use the conventions of Connes and Marcolli \cite{cm}, which distinguish between the KMS$_\infty$ states (those which are weak* limits of KMS$_\beta$ states as $\beta\to \infty$) and the ground states (those such that $z\mapsto \phi(c\sigma_z(d))$ is bounded on the upper half-plane). Neither \cite{bra-rob} nor \cite{ped} makes this distinction.

\section{A presentation}\label{sec-pres}

We describe a presentation of the Toeplitz algebra $\TT(M_L)$ like that of $\TT(\nxnx)$ in \cite[Theorem~4.1]{lr}, or, more precisely, like that of the additive boundary quotient $\TT_{\add}(\nxnx)$ in \cite[Proposition~3.3]{BaHLR}.

\begin{prop}\label{presadd}
Suppose that $A\in M_d(\Z)$ has $\det A\not= 0$, and consider the Exel system $(C(\T^d),\alpha_A,L)$. Then the Toeplitz algebra $\TT(M_L)$ is the universal $C^*$-algebra generated by a unitary representation $u:\Z^d\to U(\TT(M_L))$ and an isometry $v\in \TT(M_L)$ satisfying
\begin{itemize}
\item[(E1)] $vu_m=u_{Bm}v$, and
\smallskip
\item[(E2)] $v^*u_mv=\begin{cases} u_{B^{-1}m}&\text{if $m\in B\Z^d$}\\
0&\text{otherwise.}\end{cases}$
\end{itemize}
If $U$ is a unitary representation of $\Z^d$ in a $C^*$-algebra $C$ and $V\in C$ is an isometry satisfying \textnormal{(E1)} and \textnormal{(E2)}, then the corresponding representation $(\psi,\pi)$ of $M_L$ satisfies $U_m=\pi(\gamma_m)$ and $V=\psi(1)$.
\end{prop}

\begin{remark}
It might be helpful to see how our presentation is related to the presentation of $\TT_{\add}(\nxnx)$ in \cite{BaHLR}. The  isometries $\{v_p:p\in\primes\}$ have become the single isometry $v$, and the additive generator $s$ has been replaced by the unitary representation $u$ of the additive group $\Z^d$. The relations (T2) and (T3) in \cite[Proposition~3.3]{BaHLR} are not needed because here ``we only have one prime'' (which we will normalise to $e$ when we define our dynamics!), and the relation (Q6) in \cite[Proposition~3.3]{BaHLR} is replaced by the assumption that $u$ is a unitary representation. So we are left with (T1) and (T5), which are analogous to (E1) and (E2) respectively.

The relation (E2) implies that $\{u_mv:m\in \Sigma\}$ is a Toeplitz-Cuntz family. The analogue of the relation (Q5) used in \cite{lr} and \cite{BaHLR} is the Cuntz relation
\begin{itemize}
\item[(E3)] $1=\sum_{m\in\Sigma} (u_mv)(u_mv)^*$,
\end{itemize}
which is satisfied in the Cuntz-Pimsner algebra $\OO(M_L)$ (see Proposition~\ref{CPquotient} below).
\end{remark}

\begin{proof}[Proof of Proposition~\ref{presadd}]
The Toeplitz algebra $\TT(M_L)$ is generated by a universal representation $(i_{M_L},i_{C(\T^d)})$. It is shown in \cite[Corollary~3.3]{br} that $\TT(M_L)$ is generated by the range of $i_{C(\T^d)}$ and the single element $S:=i_{M_L}(1)$, that $(i_{C(\T^d)},S)$ is a Toeplitz-covariant representation in the sense of \cite[Definition~3.1]{br}, and  that $(\TT(M_L),i_{C(\T^d)},S)$ is universal for Toeplitz-covariant representations $(\rho, V)$ satisfying
\begin{itemize}
\item[(TC1)] $V\rho(a)=\rho(\alpha_A(a))V$, and\label{TCrels}
\smallskip
\item[(TC2)] $V^*\rho(a)V=\rho(L(a))$.
\end{itemize} 
In our system $L(1)=1$, and (TC2) implies that the operator $V$ is an isometry. 

The Stone-Weierstrass theorem implies that the characters $\gamma_m:z\mapsto z^m$ of $\T^d$ span a dense $*$-subalgebra of $C(\T^d)$, and a representation $\rho$ of $C(\T^d)$ is completely determined by the unitary representation $u:m\mapsto \rho(\gamma_m)$ of $\Z^d$. One checks that $\alpha_A(\gamma_m)=\gamma_{Bm}$, so (TC1) is equivalent to (E1). We will complete the proof by showing that (TC2) is equivalent to (E2).

To see what (TC2) says about the $u_m$, we need to compute $L(\gamma_m)$. For any $f\in C(\T^d)$ and $z=e^{2\pi ix}\in \T^d$, we can compute $L(f)(z)$ by choosing one solution $w_0$ of $\sigma_A(w)=z$, such as $w_0=e^{2\pi iA^{-1}x}$, and computing
\[
L(f)(e^{2\pi ix})=\frac{1}{|\det A|}\sum_{w\in\ker \sigma_A}f(ww_0),
\]
and hence
\[
L(\gamma_m)(e^{2\pi ix})=\frac{1}{|\det A|}\gamma_m(e^{2\pi iA^{-1}x})\sum_{w\in\ker \sigma_A}\gamma_m(w).
\]
If $\gamma_m|_{\ker\sigma_A}$ is not the identity character, then $\{\gamma_m(w):w\in\ker\sigma_A\}$ is a nontrivial subgroup of $\T$, the sum is zero, and $L(\gamma_m)=0$. So $L(\gamma_m)\not=0\Longleftrightarrow \gamma_m\in(\ker\sigma_A)^\perp$, and for such $m$, $L(\gamma_m)(e^{2\pi ix})=\gamma_m(e^{2\pi iA^{-1}x})=\gamma_{B^{-1}m}(e^{2\pi i x})$. Now
\begin{align*}
\gamma_m\in(\ker\sigma_A)^\perp
&\Longleftrightarrow e^{2\pi im^tA^{-1}n}=1\text{ for all }n\in \Z^d\\
&\Longleftrightarrow e^{2\pi i(B^{-1}m)^tn}=1\text{ for all }n\in \Z^d\\
&\Longleftrightarrow m\in B\Z^d.
\end{align*}
Thus
\[
L(\gamma_m)=\begin{cases}0&\text{unless $m\in B\Z^d$}\\
\gamma_{B^{-1}m}&\text{if $m\in B\Z^d$,}
\end{cases}
\]
and using this we can see that (TC2) is equivalent to (E2).

For the last comment, recall from \cite[\S2]{br} that the representation $(\psi,\pi)$ corresponding to the Toeplitz-covariant representation $(\rho,V)$ in the above argument is characterised by $\pi=\rho$ and $V=\psi(1)$, and that $\rho$ satisfies $\rho(\gamma_m)=u_m$. 
\end{proof}

We now want an analogous presentation of $\OO(M_L)$. To help keep things straight later, we write $\bar u:=Q\circ u$ and $\bar v:=Q(v)$.
 
\begin{prop}\label{CPquotient}
Suppose that $A\in M_d(\Z)$ has $\det A\not= 0$, and consider the Exel system $(C(\T^d),\alpha_A,L)$. Then the Cuntz-Pimsner algebra $\OO(M_L)$ is the universal $C^*$-algebra generated by a unitary representation $\bar u:\Z^d\to U(\OO(M_L))$ and an isometry $\bar v\in \OO(M_L)$ satisfying \textnormal{(E1)}, \textnormal{(E2)} and \textnormal{(E3)}.
\end{prop}

\begin{proof}
We need to prove that the unitary representation $\bar u$ and the isometry $\bar v$ satisfy (E1--3) and are universal for families satisfying these relations.
They satisfy (E1) and (E2) because $u$ and $v$ do. To see that they satisfy (E3), note that the unitary $u_m$ in Proposition~\ref{presadd} is $i_{C(\T^d)}(\gamma_m)$ and the isometry $v$ is $i_{M_L}(1)$.  We know from \cite[Lemma~2.6]{ehr} that $\{\gamma_m:m\in \Sigma\}$ is an orthonormal basis for $M_L$, so Lemma~2.5 of \cite{ehr} says that a representation $(\psi,\pi)$ is Cuntz-Pimsner covariant if and only if
\begin{equation}\label{CPcovgeneral}
1=\sum_{m\in \Sigma}\psi(\gamma_m)\psi(\gamma_m)^*.
\end{equation}
Since $\gamma_m=\gamma_m\cdot 1$ in $M_L$, we have 
\begin{equation}\label{checkCP}
Q\circ i_{M_L}(\gamma_m)=Q\circ i_{C(\T^d)}(\gamma_m)Q\circ i_{M_L}(1)=Q(u_m)Q(v), 
\end{equation} 
and Equation~\eqref{CPcovgeneral} for $(Q\circ i_{M_L}, Q\circ i_{C(\T^d)})$ reduces to (E3) for $Q\circ u$ and $Q(v)$. 

Next suppose that $U_m$ and $V$ satisfy (E1), (E2) and (E3). Then Proposition~\ref{presadd} gives a representation $(\psi,\pi)$ of $M_L$ such that $U_m=\pi(\gamma_m)$ and $V=\psi(1)$, and, in view of \eqref{checkCP},  (E3) implies that $\psi$ satisfies Equation~\eqref{CPcovgeneral}. Thus $(\psi,\pi)$ is Cuntz-Pimsner covariant, and hence factors through a representation of $\OO(M_L)$, and this representation takes $\bar u_m =Q(u_m)$ to $U_m$ and $\bar v=Q(v)$ to $V$ because $\psi\times\pi$ takes $u_m$ to $U_m$ and $v$ to $V$. 
\end{proof}

Next we want a convenient spanning family to do calculations with. Again, we are looking for something similar to what we used in \cite{lr}.

\begin{lemma} In $\TT(M_L)$ we have
\begin{equation}
(u_mv^kv^{*l}u_n^*)(u_pv^iv^{*j}u_q^*)
=\begin{cases}
u_{m+B^{k-l}(p-n)}v^{k+i-l}v^{*j}u_q^*
&\text{if $i\geq l$ and $p-n\in B^l\Z^d$}\\
u_mv^kv^{*(l+j-i)}u_{B^{j-i}(n-p)+q}^*&\text{if $i<l$ and $p-n\in B^i\Z^d$}\\
0&\text{otherwise.}
\end{cases}\label{alternatives}
\end{equation}
\end{lemma}

\begin{proof}
We compute, using first (E2) and then (E1), to get
\begin{align*}
(u_mv^kv^{*l}u_n^*)&(u_pv^iv^{*j}u_q^*)
=u_mv^kv^{*l}u_{p-n}v^iv^{*j}u_q^*\notag\\
&=\notag\begin{cases}
u_mv^ku_{B^{-l}(p-n)}v^{i-l}v^{*j}u_q^*
&\text{if $i\geq l$ and $p-n\in B^l\Z^d$}\\
u_mv^kv^{*(l-i)}u_{B^{-i}(n-p)}^*v^{*j}u_q^*&\text{if $i<l$ and $p-n\in B^i\Z^d$}\\
0&\text{otherwise}
\end{cases}\\
&=\begin{cases}
u_{m+B^{k-l}(p-n)}v^{k+i-l}v^{*j}u_q^*
&\text{if $i\geq l$ and $p-n\in B^l\Z^d$}\\
u_mv^kv^{*(l+j-i)}u_{B^{j-i}(n-p)+q}^*&\text{if $i<l$ and $p-n\in B^i\Z^d$}\\
0&\text{otherwise,}
\end{cases}
\end{align*}
as required.
\end{proof}

\begin{cor}
We have
\begin{equation}\label{spanningset}
\TT(M_L)=\clsp\big\{u_mv^kv^{*l}u_n^*:m,n\in\Z^d,\ k,l\in\N\big\}.
\end{equation}
\end{cor}

\begin{proof}
Equation~\eqref{alternatives} implies that $\newspan\{u_mv^kv^{*l}u_n^*\}$ is a $*$-algebra, and it contains all the generators of $\TT(M_L)$.
\end{proof}

\begin{remark}
Since $u$ is a unitary representation, we have $u_n^*=u_{-n}$, and the $^*$ in $u_n^*$ in \eqref{spanningset} is technically redundant. We have retained the $^*$ to emphasise the parallels between this situation and the one in \cite{lr}. It also makes formulas more symmetric, and this sometimes simplifies calculations.
\end{remark}

\section{A characterisation of KMS states}\label{sec-charKMS}

The gauge action $\gamma:\T\to \Aut \TT(M_L)$ is characterised by $\gamma_z(i_{M_L}(x))=zi_{M_L}(x)$ and $\gamma_z(i_{C(\T^d)}(f))=i_{C(\T^d)}(f)$, or equivalently by $\gamma_z(v)=zv$ and $\gamma_z(u_m)=u_m$. Our dynamics $\sigma:\R\to \Aut\TT(M_L)$ is defined in terms of the gauge action by $\sigma_t=\gamma_{e^{it}}$. Then we have
\begin{equation*}\label{sigmagen}
\sigma_t(u_mv^kv^{*l}u_n^*)=e^{it(k-l)}u_mv^kv^{*l}u_n^*,
\end{equation*} 
which since $z\mapsto e^{iz(k-l)}$ is entire implies that the spanning elements $u_mv^kv^{*l}u_n^*$ are all analytic elements. Thus a state $\phi$ on $\TT(M_L)$ is a KMS$_\beta$ state for $\sigma$ if and only if 
\begin{equation}\label{defKMSbeta}
\phi\big((u_mv^kv^{*l}u_n^*)(u_pv^iv^{*j}u_q^*)\big)
=e^{-(k-l)\beta}\phi\big((u_pv^iv^{*j}u_q^*)(u_mv^kv^{*l}u_n^*)\big)
\end{equation}
for all $m,n,p,q\in \Z^d$ and $i,j,k,l\in\N$.

The next result is an analogue of \cite[Lemma~8.3]{lr}.

\begin{prop}\label{charKMSbeta}
The system $(\TT(M_L),\sigma)$ has no KMS$_\beta$ states for $\beta< \log|\det A|$. For $\beta\geq \log|\det A|$, a state $\phi$ of $\TT(M_L)$ is a KMS$_\beta$ state if and only if 
\begin{equation}\label{charKMS}
\phi(u_mv^kv^{*l}u_n^*)=
\begin{cases}0&\text{unless $k= l$ and $m-n\in B^k\Z^d$}\\
e^{-k\beta}\phi(u_{B^{-k}(m-n)})&\text{if $k= l$ and $m-n\in B^k\Z^d$.}
\end{cases} 
\end{equation}
\end{prop}

\begin{proof}
Suppose that $\phi$ is a KMS$_\beta$ state on $(\TT(M_L),\sigma)$. Then for every $m\in\Z^d$, we have
\[
\phi(u_mvv^*u_m^*)=\phi(v^*u_m^*\sigma_{i\beta}(u_mv))=e^{-\beta}\psi(v^*u_m^*u_mv)=e^{-\beta}\psi(1)=e^{-\beta}.
\]
Since $\{u_mv:m\in\Sigma\}$ is a Toeplitz-Cuntz family with $|\det A|$ elements, we have
\[
1=\phi(1)\geq\sum_{m\in\Sigma}\phi(u_mvv^*u_m^*)=|\det A|e^{-\beta}.
\]
Thus $e^\beta\geq |\det A|$, and $\beta\geq \log|\det A|$.

Next we verify that $\phi$ satisfies \eqref{charKMS}. Applying the KMS condition twice gives
\[
\phi(u_mv^kv^{*l}u_n^*)=e^{-k\beta}\phi(v^{*l}u_n^*u_mv^k)=e^{-(k-l)\beta}\phi(u_mv^kv^{*l}u_n^*).
\]
Since $\beta\geq \log|\det A|>0$, this implies that $\phi(u_mv^kv^{*l}u_n^*)=0$ unless $k=l$. For $k=l$ we have
\begin{equation}\label{calcphi}
\phi(u_mv^kv^{*k}u_n^*)=e^{-k\beta}\phi(v^{*k}u_n^*u_mv^k)=e^{-k\beta}\phi(v^{*k}u_{m-n}v^k).
\end{equation}
From $k$ applications of (E2), we see that
\[
v^{*k}u_{m-n}v^k=\begin{cases}0&\text{unless $m-n\in B^k\Z^d$}\\
u_{B^{-k}(m-n)}&\text{if $m-n\in B^k\Z^d$,}
\end{cases}
\]
and \eqref{calcphi} implies \eqref{charKMS}.

Now we suppose that $\phi$ is a state of $\TT(M_L)$ satisfying \eqref{charKMS}, and aim to prove that $\phi$ is a KMS$_\beta$ state by verifying \eqref{defKMSbeta}. There is a certain amount of symmetry to the two nonzero alternatives in formula \eqref{alternatives}, so we may as well assume that $i\geq l$ and $p-n\in B^l\Z^d$. Then \eqref{alternatives} and \eqref{charKMS} imply that $\phi\big((u_mu^ku^{*l}u_n^*)(u_pv^iv^{*j}u_q^*)\big)$ is
\[
\begin{cases}
0&\text{unless $k+i-l=j$ and $m+B^{k-l}(p-n)-q\in B^j\Z^d$}\\
e^{-j\beta}\phi(u_{B^{-j}(m+B^{k-l}(p-n)-q)})&\text{if $k+i-l=j$ and $m+B^{k-l}(p-n)-q\in B^j\Z^d$.}
\end{cases}
\]
The right-hand side of \eqref{defKMSbeta} also vanishes unless $k+i-l=j$, so we assume this from now on. Rewriting this equation as $k-j=l-i$ shows that our assumption $i\geq l$ is equivalent to $k\leq j$. Thus when we calculate the right-hand side of \eqref{defKMSbeta}, the second alternative in \eqref{alternatives} comes into play:
\begin{align*}
\phi((&u_pv^iv^{*j}u_q^*)(u_mv^kv^{*l}u_n^*))
=\begin{cases}\phi(u_pv^iv^{*(j+l-k)}u_{B^{l-k}(q-m)+n}^*)&\text{if $q-m\in B^k\Z^d$}\\
0&\text{otherwise}\end{cases}\\
&=\begin{cases}
e^{-i\beta}\phi(u_{B^{-i}(p-B^{l-k}(q-m)-n)})&\text{if $q-m\in B^k\Z^d$ and $p-B^{l-k}(q-m)-n\in B^i\Z^d$}\\
0&\text{otherwise.}\end{cases}
\end{align*}
Multiplying this by $e^{-(k-l)\beta}$ gives the right-hand side of \eqref{defKMSbeta}, and in view of the equation $k+i-l=j$, it follows that the right-hand side of \eqref{defKMSbeta} is 
\[
\begin{cases}
e^{-j\beta}\phi(u_{B^{-i}(p-B^{l-k}(q-m)-n)})&\text{if $q-m\in B^k\Z^d$ and $p-B^{l-k}(q-m)-n\in B^i\Z^d$}\\
0&\text{otherwise.}\end{cases}
\]
Since 
\begin{align*}
B^{-i}(p-B^{l-k}(q-m)-n)&=B^{-i}(p-n)+B^{-(i+k-l)}(m-q)\\
&=B^{-j+k-l}(p-n)+B^{-j}(m-q),
\end{align*}
we have
\[
e^{-j\beta}\phi(u_{B^{-i}(p-B^{l-k}(q-m)-n)})=e^{-j\beta}\phi(u_{B^{-j}(m+B^{k-l}(p-n)-q)}),
\]
so the numbers arising on the two sides of \eqref{defKMSbeta} are the same, and it remains to check that the conditions for nonvanishing are equivalent. So we need to check that
\begin{align*}
p-n\in B^l\Z^d\text{ and }&m+B^{k-l}(p-n)-q\in B^j\Z^d\\
&\Longleftrightarrow q-m\in B^k\Z^d\text{ and }p-B^{l-k}(q-m)-n\in B^i\Z^d.
\end{align*}
Suppose that the first set of conditions holds. Then $q-m$ belongs to the coset $B^{(k-l)}(p-n)+B^j\Z^d$, which is contained in $B^k\Z^d$ because $k\leq j$ and $B^{-l}(p-n)$ is in $\Z^d$, and
\[
p-B^{l-k}(q-m)-n=B^{l-k}(B^{k-l}(p-n)-(q-m))
\]
belongs to $B^{l-k}B^j\Z^d=B^{l-k+j}\Z^d=B^i\Z^d$. So the forward implication holds, and similar arguments prove the converse. We have now proved \eqref{defKMSbeta}, and thus $\phi$ is a KMS$_\beta$ state.
\end{proof}

\section{KMS states for $\beta=\log|\det A|$.}\label{sec-Exelcp}

We begin by showing that the Exel crossed product has very few KMS states.

\begin{prop}\label{onlybeta}
If $\phi$ is a KMS$_\beta$ state on $(\OO(M_L),\sigma)$, then $\beta=\log|\det A|$.
\end{prop}

\begin{proof}
We compute using the relation (E3) and the KMS condition:
\begin{align*}
1=\phi(1)&=\phi\Big(\sum_{m\in\Sigma}\bar u_m\bar v\bar v^*\bar u_m^*\Big)=\sum_{m\in\Sigma}\phi(\bar u_m\bar v\bar v^*\bar u_m^*)\\
&=\sum_{m\in\Sigma}\phi(\bar v^*\bar u_m^*\sigma_{i\beta}(\bar u_m\bar v))
=\sum_{m\in\Sigma}e^{-\beta}\phi(\bar v^*\bar u_m^*\bar u_m\bar v)\\
&=\sum_{m\in\Sigma}e^{-\beta}\phi(1)=Ne^{-\beta},
\end{align*}
since $|\Sigma|=|\Z^d/B\Z^d|=|\det B|=|\det A|=N$.
\end{proof}

\begin{lemma}
Every KMS$_{\log N}$ state of $(\TT(M_L),\sigma)$ factors through the quotient map $Q$ of $\TT(M_L)$ onto $\OO(M_L)$. 
\end{lemma}

\begin{proof}
Suppose that $\phi$ is a KMS$_{\log N}$ state of $(\TT(M_L),\sigma)$. The formula \eqref{charKMS} implies that $\phi(u_mvv^*u_m^*)=N^{-1}$, and hence 
\[
\phi\Big(1-\sum_{m\in\Sigma}u_mvv^*u_m^*\Big)=1-NN^{-1}=0.
\]
Now the argument of \cite[Lemma~10.3]{lr} implies that $\phi$ vanishes on the ideal generated by $1-\sum_{m\in\Sigma}u_mvv^*u_m^*$. But Proposition~\ref{CPquotient} says that this ideal is the kernel of $Q$, and the result follows.
\end{proof}

The next result was first obtained by Ted Boey as an application of the general theory in~\cite{ln}. 
 
\begin{thm}\label{betalogN}
Suppose that $A\in M_d(\Z)$ has $N:=|\det A|\not=0$. Then there is a KMS$_{\log N}$ state $\phi$ of $(\OO(M_L),\sigma)$ such that 
\begin{equation}\label{KMSlogN}
\phi(\bar u_m\bar v^k\bar v^{*l}\bar u_n^*)=
\begin{cases}0&\text{unless $k=l$ and $m=n$}\\
N^{-k}&\text{if $k=l$ and $m=n$.}
\end{cases}
\end{equation} 
If $A$ is a dilation matrix, then this is the only KMS state of $(\OO(M_L),\sigma)$.
\end{thm}

We will construct the state by factoring through an expectation onto the commutative subalgebra spanned by the range projections of the generators.

\begin{lemma}\label{existexpect}
Suppose that $A\in M_d(\Z)$ has nonzero determinant. Then there is an expectation $E$ of $\OO(M_L)$ onto
\begin{equation}\label{fpa}
\OO(M_L)^\delta:=\clsp\big\{\bar u_m\bar v^k\bar v^{*k}\bar u_m^*:m\in\Z^d,\ k\in\N\big\}
\end{equation}
such that 
\begin{equation}\label{defexpect}
E(\bar u_m\bar v^k\bar v^{*l}\bar u_n^*)=
\begin{cases}0&\text{unless $k=l$ and $m=n$}\\
\bar u_m\bar v^k\bar v^{*k}\bar u_m^*&\text{if $k=l$ and $m=n$.}
\end{cases}
\end{equation} \end{lemma}

We prove this by averaging over a dual coaction, following a line of argument used in \cite{lr1} and \cite{lr} (and this explains the notation $\OO(M_L)^\delta$). We will later give a second proof which avoids the use of coactions.

\begin{proof}[First proof of Lemma~\ref{existexpect}]
The Baumslag-Solitar group $\Z[B^{-1}]\rtimes \Z$ is the semidirect product of the additive subgroup $\Z[B^{-1}]:=\bigcup_k B^{-k}\Z^d$ of $\Q$ by the action of $\Z$ by powers of $B$. We write $\epsilon$ for the canonical unitary representation of $\Z[B^{-1}]\rtimes \Z$ in $C^*(\Z[B^{-1}]\rtimes \Z)$, so that in particular
\begin{align*}
\epsilon_{(0,1)}\epsilon_{(m,0)}&=\epsilon_{(Bm,1)}=\epsilon_{(Bm,0)}\epsilon_{(0,1)}, \ \text{ and}\\
\epsilon_{(0,1)}^*\epsilon_{(m,0)}\epsilon_{(0,1)}&=\epsilon_{(B^{-1}m,0)}.
\end{align*}
These identities imply that $U_m:=\bar u_m\otimes\epsilon_{(m,0)}$ and $V:=\bar v\otimes \epsilon_{(0,1)}$ satisfy the relations (E1), (E2) and (E3), and hence give a homomorphism $\delta:=\pi_{U,V}$ of $\OO(M_L)$ into $\OO(M_L)\otimes C^*(\Z[B^{-1}]\rtimes \Z)$. One can check on generators that $\delta$ is a coaction of 
$\Z[B^{-1}]\rtimes \Z$ on $\OO(M_L)$. Since $\Z[B^{-1}]\rtimes \Z$ is amenable, averaging over this coaction gives an expectation $E$ of $\OO(M_L)$ onto the fixed-point algebra
\[
\OO(M_L)^\delta:=\big\{a\in \OO(M_L): \delta(a)=a\otimes 1=a\otimes \epsilon_{(0,0)}\big\}
\]
(see \cite[Lemma~6.5]{lr1}). Since  
\[
\delta(\bar u_m\bar v^k\bar v^{*l}\bar u_n^*)=\bar u_m\bar v^k\bar v^{*l}\bar u_n^*\otimes \epsilon_{(m,k)(n,l)^{-1}},
\]
$\bar u_m\bar v^k\bar v^{*l}\bar u_n^*$ belongs to the fixed-point algebra if and only if $(m,k)=(n,l)$, and thus $E$ satisfies \eqref{defexpect}. Since $E$ is norm-decreasing, and the $\bar u_m\bar v^k\bar v^{*l}\bar u_n^*$ span a dense subspace of $\OO(M_L)$, \eqref{defexpect} implies \eqref{fpa}.
\end{proof}

The coaction-free proof of Lemma~\ref{existexpect} involves averaging twice over actions of abelian groups. Averaging over the gauge action $\gamma:\T\to \Aut\OO(M_L)$ gives an expectation $E^\gamma$ onto the fixed-point algebra $\OO(M_L)^\gamma$; since this expectation is continuous and kills elements $\bar u_m\bar v^k\bar v^{*l}\bar u_n^*$ with $k\not= l$, we have
\begin{equation}\label{gaugefpa}
\OO(M_L)^\gamma=\clsp\big\{\bar u_m\bar v^k\bar v^{*k}\bar u_n^*:m,n\in\Z^d,\ k\in\N\big\}.
\end{equation}
For this proof of Lemma~\ref{existexpect}, we need to analyse the structure of $\OO(M_L)^\gamma$, and since we'll use this analysis elsewhere in the proof of Theorem~\ref{betalogN}, we might as well do it properly now. As a point of notation, we write
\begin{equation}\label{defSigmak}
\Sigma_k:=\big\{\mu_1+B\mu_2+\cdots B^{k-1}\mu_k:\mu\in\Sigma^k\big\},
\end{equation}
and observe that $\Sigma_k$ is a set of coset representatives for $\Z^d/B^k\Z^d$.

\begin{prop}\label{analcore} \textnormal{(a)} For each $k\geq 1$, we set
\[
C_k:=\clsp\big\{\bar u_m\bar v^k\bar v^{*k}\bar u_n^*:m,n\in\Z^d\big\}.
\]
Then the $C_k$ are $C^*$-subalgebras of $\OO(M_L)^\gamma$ satisfying $C_k\subset C_{k+1}$ and $\OO(M_L)^\gamma=\overline{\bigcup_{k=1}^\infty C_k}$.

\smallskip \textnormal{(b)} For each $k\geq 1$, $\{e_{m,n}^k:= \bar u_m\bar v^k\bar v^{*k}\bar u_n^*:m,n\in\Sigma_k\}$ is a set of nonzero matrix units which spans a matrix algebra $M_{\Sigma_k}(\C)$. The respresentation $\bar u$ of $\Z^d$ in $\OO(M_L)$ maps $B^k\Z^d$ into $C_k$, and every $\bar u_{B^km}$ belongs to the commutant of $M_{\Sigma_k}(\C)$ in $C_k$.

\smallskip \textnormal{(c)} The inclusion $\iota_k$ of $M_{\Sigma_k}(\C)$ in $\OO(M_L)^\gamma$ and the integrated form $\pi_{\bar u,k}$ of $\bar u|_{B^k\Z^d}$ give an isomorphism $\iota_k\otimes \pi_{\bar u,k}$ of $M_{\Sigma_k}(\C)\otimes C^*(B^k\Z^d)$ onto $C_k$ which carries $e_{m,n}^k\otimes \epsilon_{B^kp}$ into $\bar u_{m+B^kp}\bar v^k\bar v^{*k}\bar u_{n}^*$.
\end{prop} 

\begin{proof}
Calculations using the relations (E1) and (E2) show that
\begin{equation}\label{calcinCk}
(\bar u_m\bar v^k\bar v^{*k}\bar u_n^*)(\bar u_p\bar v^k\bar v^{*k}\bar u_q^*)
=\begin{cases}
0&\text{unless $p-n\in B^k\Z^d$}\\
\bar u_{m+p-n}\bar v^k\bar v^{*k}\bar u_q^*&\text{if $p-n\in B^k\Z^d$,}
\end{cases}
\end{equation}
which implies that $C_k$ is a $C^*$-subalgebra. The Cuntz relation (E3) implies that
\begin{align*}
\bar u_m\bar v^k\bar v^{*k}\bar u_n^*
&=\sum_{p\in\Sigma}\bar u_m\bar v^k(\bar u_p\bar v\bar v^{*}\bar u_p^*)\bar v^{*k}\bar u_n^*=\sum_{p\in\Sigma}\bar u_{m+B^kp}\bar v^{k+1}\bar v^{*(k+1)}\bar u_{n+B^kp}^*,
\end{align*}
which in view of the definition of $\Sigma_{k+1}$ implies that $C_k\subset C_{k+1}$. Now $\OO(M_L)^\gamma=\overline{\bigcup_{k=1}^\infty C_k}$ follows from \eqref{gaugefpa}, and we have proved (a).

Since $\{S_m:=\bar u_m\bar v:m\in\Sigma\}$ is a Cuntz family, for each fixed $k$ the products $\{S_\mu=S_{\mu_1}\cdots S_{\mu_k}:\mu\in\Sigma^k\}$ form a Cuntz family; since
\[
S_\mu=(\bar u_{\mu_1}\bar v)(\bar u_{\mu_2}\bar v)\cdots(\bar u_{\mu_k}\bar v)=\bar u_{\mu_1+B\mu_2+\cdots B^{k-1}\mu_k}\bar v^k,
\]
this Cuntz family is precisely $\{\bar u_m\bar v^k:m\in \Sigma_k\}$, and it follows that the $e_{m,n}^k$ are nonzero matrix units. The relation (E1) implies that $\bar u_{B^kp}$ commutes with every $e_{m,n}^k$, which gives (b).

Since the representation of $C^*(B^k\Z^d)$ in $C(\T^d)$ is faithful, and since $j_{C(\T^d)}:C(\T^d)\to C(\T^d)\rtimes_{\alpha_A,L}\N=\OO(M_L)$ is injective (by \cite[Corollary~4.3]{br}, for example), $\pi_{\bar u,k}$ is injective, and hence so is the representation $\iota_k\otimes \pi_{\bar u,k}$ of $M_{\Sigma_k(\C)}\otimes C^*(B^k\Z^d)=M_{\Sigma_k}(C^*(B^k\Z^d))$. It is surjective because every $m\in\Z^d$ can be written uniquely as $m'+B^km''$ for some $m'\in\Sigma_k$, and then
\[
\bar u_m\bar v^k\bar v^{*k}\bar u_n^*=\bar u_{B^km''}(\bar u_{m'}\bar v^k\bar v^{*k}\bar u_{n'}^*)\bar u_{B^kn''}^*
\]
is in $C_k$ because each matrix unit and each $\bar u_{B^km''}$ are. To see the last assertion about what $\iota_k\otimes \pi_{\bar u,k}$ does to $e_{m,n}^k\otimes \epsilon_{B^kp}$, recall that the representation $\phi\otimes\psi$ of a tensor product $C\otimes D$ coming from commuting representations $\phi$ of $C$ and $\psi$ of $D$ takes $c\otimes d$ to $\phi(c)\psi(d)=\psi(d)\phi(c)$ (see \cite[Theorem~B.2]{tfb}, for example).
\end{proof}

\begin{cor}\label{2ndaction}
There is a continuous action $\tau$ of $\T^d$ on $\OO(M_L)^\gamma$ such that\begin{equation}\label{deftau}
\tau_z(\bar u_m\bar v^k\bar v^{*k}\bar u_n^*)=z^{m-n}\bar u_m\bar v^k\bar v^{*k}\bar u_n^*\ \text{ for $m,n\in\Z^d$, $k\in\N$.}
\end {equation}
\end{cor}

\begin{proof}
There is a continuous action $\eta$ of $\T^d$ on $M_{\Sigma_k}(\C)$ such that $\eta_z(e_{m,n}^k)=z^{m-n}e_{m,n}^k$ --- indeed, $\eta_z$ is conjugation by the unitary $\sum_{m\in\Sigma_k}z^me_{m,m}^k$. Lifting the dual action of $(B^k\Z^d)^\wedge=\T^d/(B^k\Z^d)^\perp$ to $\T^d$ gives an action $\zeta$ of $\T^d$ on $C^*(B^k\Z^d)$ which multiplies the generator $\epsilon_{B^km}$ by $z^{B^km}$. Pulling the action $\eta\otimes \zeta$ on $M_{\Sigma_k(\C)}\otimes C^*(B^k\Z^d)$ over to $\OO(M_L)^\gamma$ under the isomorphism of Proposition~\ref{analcore}(c) gives an action $\tau^k$ on $C_k$ which satisfies \eqref{deftau} (for fixed $k$). A calculation using the Cuntz relation (E3) shows that the automorphisms $\tau^k_z$ combine to give an automorphism $\tau_z$ of $\bigcup_{k=1}^\infty C_k$, which is isometric because each $\tau^k_z$ is, and hence extends to an automorphism of $\OO(M_L)^\gamma$. Continuity follows from the continuity of scalar multiplication.
\end{proof}

\begin{proof}[Second proof of Lemma~\ref{existexpect}]
Averaging over the action $\tau$ of Corollary~\ref{2ndaction} gives an expectation $E^\tau$ of $\OO(M_L)^\gamma$ onto $\OO(M_L)^\delta$, and $E:=E^\tau\circ E^\gamma$ has the required properties.
\end{proof}

\begin{proof}[Proof of existence in Theorem~\ref{betalogN}]
The description of $\OO(M_L)^\delta$ in \eqref{fpa} shows that each $\sigma_t$ is the identity on $\OO(M_L)^\delta$, so it suffices to find a trace $\tau$ on $\OO(M_L)^\delta$ such that $\tau(\bar u_m\bar v^k\bar v^{*k}\bar u_m^*)=N^{-k}$, and then Proposition~\ref{charKMSbeta} implies that $\tau\circ E$ is a KMS$_{\log N}$ state on $(\OO(M_L),\sigma)$. 

We can write each $m\in \Z^d$ uniquely as $m'+B^km''$ for some $m'\in \Sigma_k$, and then part~(b) of Proposition~\ref{analcore} implies that
\begin{align*}
\bar u_m\bar v^k\bar v^{*k}\bar u_m^*&=\bar u_{B^km''}(\bar u_{m'}\bar v^k\bar v^{*k}\bar u_{m'}^*)\bar u_{B^km''}^*\\
&=\bar u_{B^km''}\bar u_{B^km''}^*(\bar u_{m'}\bar v^k\bar v^{*k}\bar u_{m'}^*)\\
&=\bar u_{m'}\bar v^k\bar v^{*k}\bar u_{m'}^*.
\end{align*}
Now part~(a) of Proposition~\ref{analcore} implies that
\[
D_k:=\newspan\big\{\bar u_m\bar v^k\bar v^{*k}\bar u_m^*:m\in\Z^d\big\}=\newspan\big\{\bar u_{m'}\bar v^k\bar v^{*k}\bar u_{m'}^*:m'\in\Sigma_k\big\}
\]
is a finite-dimensional commutative $C^*$-algebra, and that $D_k$ has a normalised trace $\tau_k$ satisfying $\tau_k(\bar u_m\bar v^k\bar v^{*k}\bar u_m^*)=N^{-k}$. The Cuntz relation (E3) implies that $D_k\subset D_{k+1}$, and the normalised traces $\tau_k$ combine to give a trace $\tau$ on $\OO(M_L)^\delta=\overline{\bigcup_{k}D_k}$ such that $\tau(\bar u_m\bar v^k\bar v^{*k}\bar u_m^*)=N^{-k}$. Then, as foreshadowed above, $\phi:=\tau\circ E$ is a KMS$_{\log N}$ state on $(\OO(M_L),\sigma)$ satisfying \eqref{KMSlogN}.
\end{proof}

For the proof of uniqueness, we need a standard fact about dilation matrices.

\begin{lemma}\label{dilpure}
If $B$ is an integer dilation matrix, then $\bigcap_{k=1}^\infty B^k\Z^d=\{0\}$.
\end{lemma}

\begin{proof} Suppose that $m\in \bigcap_{k=0}^\infty B^k\Z^d$. Then $B^{-k}m$ belongs to $\Z^d$ for every $k$, and since we know from \cite[Lemma~4.12]{ehr}, for example, that  $\|B^{-k}m\|\to 0$ as $k\to\infty$, we must have $B^km=0$ for large $k$, and $m=0$.
\end{proof}

\begin{proof}[Proof of uniqueness in Theorem~\ref{betalogN}]
Suppose that $\phi$ is a KMS state of $(\OO(M_L),\sigma)$. Proposition~\ref{onlybeta} implies that $\phi$ has inverse temperature $\beta=\log N$. We need to prove that $\phi$ satisfies \eqref{KMSlogN}, and comparing \eqref{KMSlogN} with \eqref{charKMS} (which we know holds with $e^{-k\beta}=N^{-k}$) shows that we need to prove that $\phi(\bar u_n)=0$ for all nonzero $n$. So suppose $n\in\Z^d$ and $n\not= 0$. Lemma~\ref{dilpure} implies that there is a smallest integer $k$ such that $n$ does not belong to $B^{k}\Z^d$.  Then, recalling from the proof of existence that $\{\bar u_m\bar v^k:m\in \Sigma_k\}$ is a Cuntz family in $\OO(M_L)$, we have
\[
\phi(\bar u_n)=\phi\Big(\bar u_n\sum_{m\in\Sigma_k}\bar u_m\bar v^k\bar v^{*k}\bar u_m^*\Big)=\sum_{m\in\Sigma_k}\phi(\bar u_{n+m}\bar v^k\bar v^{*k}\bar u_m^*),
\]
which vanishes by \eqref{charKMS} because $(n+m)-m=n$ is not in $B^k\Z^d$ for every $m$.
\end{proof}

It seems to be quite easy to find representations of $\OO(M_L)$, and we describe an interesting one in the following example (which was one of our reasons for becoming interested in the $C^*$-algebras associated to dilation matrices in the first place \cite{ehr}). However, it does not seem to be so easy to find natural Hilbert space representations of $\OO(M_L)$ in which the KMS$_{\log|\det A|}$ state is a vector state.

\begin{example}\label{BrattJorgrep} The operators $V$ and $U_m$ on $L^2(\T^d)$ defined by 
\[
(V\xi)(z)=\xi(\sigma_A(z))\ \text{ and }\ (U_m\xi)(z)=z^m\xi(z)
\]
satisfy (E1--3), and hence give a representation of $\OO(M_L)$ on $L^2(\T^d)$. The Cuntz family $\{U_mV:m\in \Sigma\}$ is one of the sort studied by Bratteli and Jorgensen in the context of wavelets \cite{bj}, or more precisely, one of the more general sort studied in~\cite{aijln}. 

To make the connection, note that the characters $\{\gamma_n:n\in \Sigma\}$ form an orthonormal basis for the right Hilbert module $M_L$, or what is called in \cite{aijln} a ``filter bank for dilation by $A$''.  It is shown in \cite[Proposition~2.2]{aijln} that any filter bank $\{m_i:0\leq i<N\}$ gives rise to a Cuntz family $S_i:=M(m_i)V$, where $M(f)$ is the operator of multiplication by $f\in C(\T^d)$. In the construction of wavelets, the more interesting filter banks are those in which $m_0$ is ``low-pass'', which implies in particular that $m_0(1)=N^{1/2}$ and $m_i(1)=0$ for $i>0$ (see \cite[Example~4.2]{aijln}); the filters $\gamma_n$ satisfy $|\gamma_n|\equiv1$, and hence are ``all-pass''. 
\end{example}

\begin{remark} Exel has previously studied KMS states on Exel crossed products \cite{e2}, and we now reconcile our result with his \cite[Proposition~9.2]{e2}. The situation in \cite{e2} is more general than ours, but applies with $h=e1$ and $E=\alpha\circ L$, which is easily seen to be an expectation of $C(\T^d)$ onto the range of $\alpha$; since our orthonormal basis for $M_L$ is a quasi-basis, $E$ has finite type with index $N:=|\det A|$ (strictly speaking, $\ind E$ is the element $N1$ of $C(\T^d)$). Exel proved in \cite[Theorem~8.9]{e2} that there is an expectation $G:C(\T^d)\rtimes_{\alpha,L}\N\to j_{C(\T^d)}(C(\T^d))$ such that $G(\bar u_m\bar v^k\bar v^{*l}\bar u^*_n)=\delta_{k,l}N^{-k}\bar u_m\bar u^*_n$; in our situation, it is quite easy to check directly that $G$ is given by first averaging over the gauge action $\gamma$, and then combining the expectations $G_k$ on $C_k:=\clsp\{\bar u_m\bar v^k\bar v^{*k}\bar u_n^*\}$ defined by $G_k(T)=N^{-k}\sum_{p\in\Sigma_k}\bar u_pT\bar u_p^*$ to get $G$ on $(C(\T^d)\rtimes\N)^\gamma=\overline{\bigcup_{k\geq 0}C_k}$ (see \cite[Corollary~7.5]{ev}). Then \cite[Proposition~9.2]{e2} implies that the KMS$_\beta$ states on $C(\T^d)\rtimes_{\alpha,L}\N$ have the form $\phi\circ G$, where $\phi$ is a trace on $C(\T^d)$ satisfying $\phi(f)=e^{-\beta}N\phi(L(f))$ for $f\in C(\T^d)$.

Traces on $C(\T^d)$ are given by measures $\mu$, and Exel's condition says that $\mu$ satisfies
\begin{equation}\label{Exelschar}
\int f\,d\mu=e^{-\beta}\int_{\T^d}\sum_{\sigma_A(w)=z}f(w)\,d\mu(z)\ \text{ for $f\in C(\T^d)$.}
\end{equation}
It follows from \cite[Lemma~2.3]{aijln}, for example, that the Haar measure $\lambda$ on $\T^d$ satisfies \eqref{Exelschar} with $1=e^{-\beta}N$, and since $\bar u_m\in C(\T^d)\rtimes\N$ is the image of the function $z^n$ in $C(\T^d)$, the corresponding KMS$_{\log N}$ state $\psi$ on $C(\T^d)\rtimes_{\alpha_A,L}\N$ satisfies
\begin{align*}
\psi(\bar u_m\bar v^k\bar v^{*l}\bar u_n^*)
&=
\begin{cases}0&\text{unless $k=l$}\\
\int_{\T^d} N^{-k}z^mz^{-n}\,d\lambda(z)&\text{if $k=l$}
\end{cases}\\
&=
\begin{cases}0&\text{unless $k=l$ and $m=n$}\\
N^{-k}&\text{if $k=l$ and $m=n$.}
\end{cases}
\end{align*}
Thus Exel's result also gives the KMS$_{\log N}$ state described in Theorem~\ref{betalogN}, even though his state was obtained by factoring through a different expectation on $C(\T^d)\rtimes_{\alpha_A,L}\N$.
\end{remark}

\section{Existence of KMS states for $\beta>\log|\det A|$.}\label{sec-existence}

Our goal here is to prove the existence of KMS$_\beta$ states for $\beta>\log|\det A|$.  Note that, when $A$ is a dilation matrix, Lemma~\ref{dilpure} implies that the sum on the right-hand side of \eqref{formpsibetamu} is finite.

\begin{prop}\label{existforbeta>1}
Suppose that $A\in M_d(\Z)$ satisfies $\det A\not= 0$ and that $\beta> \log|\det A|$. Then for each probability measure $\mu$ on $\T^d$, there is a KMS$_\beta$ state $\psi=\psi_{\beta,\mu}$ of $(\TT(M_L),\sigma)$ such that $\psi(u_mv^kv^{*l}u_n^*)$ vanishes unless $k= l$ and $m-n\in B^k\Z^d$, and
\begin{equation}\label{formpsibetamu}
\psi(u_mv^kv^{*l}u_n^*)=
(1-|\det A|e^{-\beta})\sum_{\{j\geq k\;:\;m-n\in B^j\Z^d\}}|\det A|^{j-k}e^{-j\beta}\int_{\T^d} z^{B^{-j}(m-n)}\,d\mu(z)
\end{equation}
when $k= l$ and $m-n\in B^k\Z^d$.
\end{prop}

We use the representation $M$ of $C(\T^d)$ by multiplication operators on $L^2(\T^d,d\mu)$, and use the same notation for the corresponding unitary representation of $\Z^d$, so that $M_m:=M(\gamma_m)$. For each $j\in \N$, we have a unitary representation $M\circ B^{-j}$ of the subgroup $B^j\Z^d$ of $\Z^d$, and we denote by $\HH_j$ the Hilbert space of the induced representation $\Ind_{B^j\Z^d}^{\Z^d}M\circ B^{-j}$. Our state $\psi_{\beta,\mu}$ will be built from vector states for a representation $\pi_\mu$ of $\TT(M_L)$ on $\HH_\mu:=\bigoplus_{j=0}^\infty\HH_j$.

We will need to do some calculations in the Hilbert spaces $\HH_j$, and for this it is convenient to use the sets $\Sigma_j$ described in \eqref{defSigmak}; for $g\in \Z^d/B^j\Z^d$, we write $c_j(g)$ for the element of $\Sigma_j$ such that $c_j(g)\in g$. Then (from \cite[page~296]{tfb}, for example) $\HH_j$ is the completion of the space 
\[
\VV_c:=\big\{\xi:\Z^d\to L^2(\T^d,d\mu)\text{ such that } \xi(m-n)=M_{B^{-j}n}(\xi(m)) \text{ for $n\in B^j\Z^d$}\big\}
\]
in the inner product defined by
\[
(\xi\,|\,\eta)=\sum_{g\in\Z^d/B^j\Z^d}\big(\xi(c_j(g))\,|\,\eta(c_j(g))\big)=\sum_{g\in\Z^d/B^j\Z^d}\int_{\T^d}\xi(c_j(g))(z)\overline{\eta(c_j(g))(z)}\,d\mu(z).
\]
(Although we have used the cross-section $c_j$ to get a useful formula for the inner product, the translation condition on $\xi$ and $\eta$ means that this inner product does not depend on the choice of $c_j$.) Then the induced representation acts on $\HH_j$ by 
\[
\big((\Ind_{B^j\Z^d}^{\Z^d}M\circ B^{-j})_m\xi\big)(n)=\xi(n-m).
\]

We now take $U$ to be the unitary representation of $\Z^d$ on $\HH_\mu$ defined by
\[
U:={\textstyle \bigoplus_{j=0}^\infty \big(\Ind_{B^j\Z^d}^{\Z^d}M\circ B^{-j}\big).}
\]
For each $j\geq 0$ and $\xi\in\HH_j$, we define
\[
(V_j\xi)(m)=\begin{cases}0&\text{unless $m\in B\Z^d$}\\
\xi(B^{-1}m)&\text{if $m\in B\Z^d$;}\end{cases}
\]
a quick calculation shows that $V_j\xi$ belongs to $\HH_{j+1}$. The $V_j$ combine to give an isometry $V$ on $\HH_\mu=\bigoplus_j\HH_j$, and the adjoint $V^*$ is given on $\HH_{j+1}$ by the formula $(V^*\xi)(n)=\xi(Bn)$. Calculations show that the pair $(U,V)$ satisfies (E1) and (E2), and hence there is a representation $\pi_\mu$ of $\TT(M_L)$ on $\HH_\mu$ such that $\pi_\mu(u_m)=U_m$ and $\pi_\mu(v)=V$.

We now let $e_{0,0}$ be the constant function $1$ viewed as a unit vector in $\HH_0=L^2(\T^d,d\mu)$. For $j\in \N$ and $g\in \Z^d/B^j\Z^d$, we define $e_{j,g}:=U_{c_j(g)}V^je_{0,0}$, so that for each $j$, 
\[
\big\{e_{j,g}:g\in \Z^d/B^j\Z^d\big\}
\]
is an orthonormal set of $|\det B|^j=|\det A|^j$ vectors in $\HH_j$. We view them as elements of $\HH_\mu$ by adding $0$s in the other summands. Inspired by the proof of \cite[Proposition~9.3]{lr}, we define
\[
\psi(T):=(1-|\det A|e^{-\beta})\sum_{j=0}^\infty\sum_{g\in\Z^d/B^j\Z^d} e^{-j\beta}\big(\pi_\mu(T)e_{j,g}\,|\,e_{j,g}\big).
\]
Summing the geometric series $\sum_j(|\det A|e^{-\beta})^j$ shows that this series converges in norm in $\TT(M_L)^*$, and that the sum is a state $\psi$ of $\TT(M_L)$.

Next we fix $m,n\in \Z^d$ and $k,l\in \N$, and verify the formula for $\psi(u_mv^kv^{*l}u_n^*)$. Then
\[
V^{*l}U_n^*e_{j,g}=V^{*l}U_n^*U_{c_j(g)}V^je_{0,0}
=\begin{cases}0&\text{unless $l\leq j$}\\
V^{*l}U_{c_j(g)-n}V^je_{0,0}&\text{if $l\leq j$}
\end{cases}
\]
belongs to $\HH_{j-l}$, and hence
\begin{align*}
\big(\pi_\mu(u_mv^kv^{*l}u_n^*)&e_{j,g}\,|\,e_{j,g}\big)
=\big(V^{*l}U_n^*e_{j,g}\,|\,V^{*k}U_m^*e_{j,g}\big)\\
&=\begin{cases}
0&\text{unless $k=l\leq j$}\\
\big(V^{*k}U_{c_j(g)-n}V^je_{0,0}\,|\,V^{*k}U_{c_j(g)-m}V^je_{0,0}\big)
&\text{if $k=l\leq j$.}
\end{cases}
\end{align*}
We now recall that $\HH_{j-k}$ is the Hilbert space of the representation $\Ind_{B^{j-k}\Z^d}^{\Z^d}(M\circ B^{j-k})$, and hence
\begin{align}
\big(V^{*k}&U_{c_j(g)-n}V^je_{0,0}\,|\,V^{*k}U_{c_j(g)-m}V^je_{0,0}\big)\notag\\
&=\sum_{h\in\Z^d/B^{j-k}\Z^d}\big(V^{*k}U_{c_j(g)-n}V^je_{0,0}(c_{j-k}(h))\,|\,V^{*k}U_{c_j(g)-m}V^je_{0,0}(c_{j-k}(h))\big)\notag\\
&=\sum_{h\in\Z^d/B^{j-k}\Z^d}\big(U_{c_j(g)-n}V^je_{0,0}(B^kc_{j-k}(h))\,|\,U_{c_j(g)-m}V^je_{0,0}(B^kc_{j-k}(h))\big)\notag\\
&=\sum_{h\in\Z^d/B^{j-k}\Z^d}\big(V^je_{0,0}(B^kc_{j-k}(h)-c_j(g)+n)\,|\,V^je_{0,0}(B^kc_{j-k}(h)-c_j(g)+m)\big).\label{formforip2}
\end{align}
The $h$-summand vanishes unless both
\begin{equation}\label{condsongh}
B^kc_{j-k}(h)-c_j(g)+n\in B^j\Z^d\ \text{ and }\ B^kc_{j-k}(h)-c_j(g)+m\in B^j\Z^d.
\end{equation}
As a function in the Hilbert space 
\[
\HH_0=\HH(\Ind_{\Z^d}^{\Z^d}M)=\big\{\xi:\Z^d\to L^2(\T^d,d\mu)\text{ such that }\xi(-n)=M_n\xi(0)\big\},
\]
$e_{0,0}$ satisfies $e_{0,0}(q)(z)=z^{-q}$, and $(V^je_{0,0})(B^jq)(z)=e_{0,0}(q)(z)=z^{-q}$. Thus, when both criteria in \eqref{condsongh} are satisfied, we have
\begin{align}
\big(\pi_\mu(u_mv^kv^{*k}u_n^*)e_{j,g}\,&|\,e_{j,g}\big)=\big(V^{*k}U_{c_j(g)-n}V^je_{0,0}\,|\,V^{*k}U_{c_j(g)-m}V^je_{0,0}\big)\label{calcip1}\\
&=\int_{\T^d} z^{-B^{-j}(B^kc_{j-k}(h)-c_j(g)+n)}\overline{z^{-B^{-j}(B^kc_{j-k}(h)-c_j(g)+m)}}\,d\mu(z)\notag\\
&=\int_{\T^d}z^{B^{-j}(m-n)}\,d\mu(z).\label{calcip}
\end{align}
(Notice that when \eqref{condsongh} holds, we have $m-n\in B^j\Z^d$, so the last integral makes sense.) For each pair $m,n$ such that $m-n$ is in $B^j\Z^d$, and each $h$ in $\Z^d/B^{j-k}\Z^d$, there is exactly one $g$ such that \eqref{condsongh} holds. Thus, using \eqref{formforip2} to view
\begin{equation}\label{collapsiblesum}
\sum_{g\in \Z^d/B^j\Z^d}e^{-j\beta}\big(\pi_\mu(u_mv^kv^{*l}u_n^*)e_{j,g}\,|\,e_{j,g}\big)
\end{equation}
as a sum over $g\in \Z^d/B^j\Z^d$ and $h\in \Z^d/B^{j-k}\Z^d$, we find that \eqref{collapsiblesum} has exactly 
\[
|\Z^d/B^{j-k}\Z^d|= |\det A|^{j-k}
\] 
nonzero terms, each of which is equal to \eqref{calcip}. Thus $\psi(u_mv^kv^{*l}u_n^*)$ vanishes unless $k=l$ and $m-n\in B^k\Z^d$, and then equals
\begin{equation}\label{ansforpsi}
(1-|\det A|e^{-\beta})\sum_{\{j\geq k\;:\;m-n\in B^j\Z^d\}}|\det A|^{j-k}e^{-j\beta}\int_{\T^d}z^{B^{-j}(m-n)}\,d\mu(z),
\end{equation}
as stated in the Proposition.

We still need to prove that $\psi$ is a KMS$_\beta$ state, and we will do this using Proposition~\ref{charKMSbeta}. So we need to compute $e^{-k\beta}\psi(u_{B^{-k}(m-n)})$ under the assumption that $m-n\in B^k\Z^d$. We have already done most of the work: the calculation \eqref{calcip1} shows that
\begin{align*}
e^{-k\beta}&\psi(u_{B^{-k}(m-n)})=e^{-k\beta}(1-|\det A|e^{-\beta})\sum_{j'=0}^\infty\sum_{g\in \Z^d/B^{j'}\Z^d}e^{-j'\beta}\big(U_{B^{-k}(m-n)}e_{j',g}\,|\,e_{j',g}\big)\\
&=e^{-k\beta}(1-|\det A|e^{-\beta})\sum_{\{j'\;:\;B^{-k}(m-n)\in B^{j'}\Z^d\}}\ \sum_{g\in \Z^d/B^{j'}\Z^d}e^{-j'\beta}\int_{\T^d}z^{B^{-j'}B^{-k}(m-n)}\,d\mu(z)\\
&=e^{-k\beta}(1-|\det A|e^{-\beta})\sum_{\{j'\;:\;B^{-k}(m-n)\in B^{j'}\Z^d\}}|\det A|^{j'}e^{-j'\beta}\int_{\T^d}z^{B^{-j'}B^{-k}(m-n)}\,d\mu(z),
\end{align*}
which reduces to \eqref{ansforpsi} on writing $j=j'+k$. Thus Proposition~\ref{charKMSbeta} implies that $\psi$ is a KMS$_\beta$ state, and this completes the proof of Proposition~\ref{existforbeta>1}.

\section{Parametrisation of KMS$_\beta$ states}\label{sec-par}

\begin{prop}\label{parametrisation}
Suppose that $A\in M_d(\Z)$ has nonzero determinant and $\beta>\log|\det A|$. Then the map $\mu\mapsto \psi_{\beta,\mu}$ of Proposition~\ref{existforbeta>1} is an affine homeomorphism of the simplex $P(\T^d)$ of probability measures onto the simplex of KMS$_\beta$ states for $(\TT(M_L),\sigma)$. 
\end{prop}

As in \cite[\S10]{lr}, the crux of the argument is a reconstruction formula which allows us to recover a KMS$_\beta$ state from its ``conditioning'' $\phi_P$ to a corner $P\TT(M_L)P$. In the present situation, though, the projection
\[
P:=1-\sum_{g\in \Z^d/B\Z^d} u_{c(g)}vv^*u_{c(g)}^*=\prod_{g\in \Z^d/B\Z^d}(1-u_{c(g)}vv^*u_{c(g)}^*)
\]
belongs to $\TT(M_L)$, so we don't need to resort to spatial arguments to make sense of the conditioning: we can just define
\[
\phi_P(a)=\frac{1}{1-|\det A|e^{-\beta}}\phi(PaP), 
\]
and then since $\phi(u_{c(g)}vv^*u_{c(g)}^*)=e^{-\beta}$, the normalising factor ensures that $\phi_P$ is a state of $\TT(M_L)$. We can now state our reconstruction formula.

\begin{prop}\label{reconstruct}
Suppose that $\beta>\log|\det A|$, and that $\phi$ is a KMS$_\beta$ state on $\TT(M_L)$. Then for every $a\in \TT(M_L)$ we have
\begin{equation}\label{reconform}
\phi(a)=\lim_{n\to\infty}(1-|\det A|e^{-\beta})\sum_{j=0}^n\sum_{g\in \Z^d/B^j\Z^d} e^{-j\beta}\phi_P\big(v^{*j}u_{c_j(g)}^*au_{c_j(g)}v^j\big).
\end{equation}
\end{prop}

Convergence of the limit in Proposition~\ref{reconstruct} will follow from the following simple lemma:

\begin{lemma}\label{approxstate}
Suppose that $\phi$ is a state of a unital $C^*$-algebra $A$, and that $\{p_n\}$ is a sequence of projections in $A$ such that $\phi(p_n)\to 1$. Then $\phi(p_nap_n)\to \phi(a)$ for every $a\in A$.
\end{lemma}

\begin{proof}
We know that $\phi(1-p_n)=1-\phi(p_n)\to 0$, so the Cauchy-Schwarz inequality for $\phi$ implies that $\phi(a(1-p_n))\to 0$ for all $a\in A$. Another application of the Cauchy-Schwarz inequality shows that $\phi(p_na(1-p_n))\to 0$ also, so
\[
\phi(a)-\phi(p_nap_n)=\phi(a(1-p_n))+\phi((1-p_n)ap_n)\to 0.\qedhere
\]
\end{proof}

When we apply Lemma~\ref{approxstate}, the projections $p_n$ will be sums of the projections in the next proposition.

\begin{prop}\label{Psorthog} For $j\in \N$ and $g\in \Z^d/B^p\Z^d$ we define 
\[
P_{j,g}:=u_{c_j(g)}v^jPv^{*j}u_{c_j(g)}^*.
\]
Then the $P_{j,g}$ are mutually orthogonal projections in $\TT(M_L)$. 
\end{prop}

The proposition follows from the next lemma.

\begin{lemma}\label{calcorthog}
For each pair $(j,g)$ and $(l,h)$ we have
\[
Pv^{*j}u_{c_j(g)}^*u_{c_j(h)}v^lP=
\begin{cases}P&\text{if $j=l$ and $g=h$}\\
0&\text{otherwise.}
\end{cases}\]
\end{lemma}

\begin{proof}
If $j\not=l$, say $j<l$, then for $m\in \Z^d$ we have
\[
Pv^{*j}u_mv^lP=\begin{cases}
Pu_{B^{-j}m}v^{l-j}P&\text{if $m\in B^j\Z^d$}\\
0&\text{otherwise.}
\end{cases}
\]
Now every $n\in \Z^d$ (including $n=B^{-j}m$) has the form $n=c(n)+Bk$, so
\[
Pu_nv^{l-j}P=Pu_{c(n)+Bk}vv^{l-j-1}P=Pu_{c(n)}vu_kv^{l-j-1}P,
\]
which vanishes because $P$ contains the factor $(1-u_{c(n)}vv^*u_{c(n)}^*)$. So $Pv^{*j}u_mv^lP$ vanishes when $j\not=l$, and for $j=l$, $Pv^{*j}u_{c_j(h)-c_j(g)}v^jP$ vanishes unless $c_j(h)-c_j(g)$ belongs to $B^j\Z^d$, which occurs precisely when $g=h$ in $\Z^d/B^j\Z^d$.
\end{proof}

\begin{proof}[Proof of Proposition~\ref{reconstruct}]
We aim to apply Lemma~\ref{approxstate} with 
\[
p_n:=\sum_{j=0}^n\sum_{g\in \Z^d/B^j\Z^d} P_{j,g},
\]
which is a projection by Proposition~\ref{Psorthog}. So we need to compute $\phi(p_n)$, which we do using the KMS condition:
\begin{align*}
\phi(p_n)&=\sum_{j=0}^n\sum_{g\in \Z^d/B^j\Z^d} \phi(P_{j,g})=\sum_{j=0}^n\sum_{g\in \Z^d/B^j\Z^d} \phi\big(u_{c_j(g)}v^jPv^{*j}u_{c_j(g)}^*\big)\\
&=\sum_{j=0}^n\sum_{g\in \Z^d/B^j\Z^d} e^{-j\beta}\phi\big(Pv^{*j}u_{c_j(g)}^*u_{c_j(g)}v^jP\big)\\
&=\phi(P)\sum_{j=0}^n |\det A|^je^{-j\beta}\qquad\text{(by Lemma~\ref{calcorthog})}\\
&=(1-|\det A|e^{-\beta})\sum_{j=0}^n |\det A|^je^{-j\beta},
\end{align*}
which on summing the geometric series converges to $1$ as $n\to \infty$. So Lemma~\ref{approxstate} implies that for every $a\in \TT(M_L)$, we have
\[
\phi(a)=\lim_{n\to\infty} \sum_{j,l=0}^n\sum_{g\in \Z^d/B^j\Z^d}\sum_{h\in \Z^d/B^l\Z^d}\phi(P_{j,g}aP_{l,h}).
\]
Applying the KMS condition shows that this sum is
\begin{align*}
\lim_{n\to\infty}\sum_{j,l=0}^n\sum_{g\in \Z^d/B^j\Z^d}\sum_{h\in \Z^d/B^l\Z^d}e^{-j\beta}\phi\big(Pv^{*j}u_{c_j(g)}^*au_{c_l(h)}v^lPv^{*l}u_{c_l(h)}^*(u_{c_j(g)}v^jP)\big),
\end{align*}
and it follows from Lemma~\ref{calcorthog} that the summands are zero unless $j=l$ and $g=h$, in which case the right-hand factor $Pv^{*l}u_{c_l(h)}^*u_{c_j(g)}v^jP$ collapses to $P$, and we recover the desired formula \eqref{reconform}.
\end{proof}

\begin{proof}[Proof of Proposition~\ref{parametrisation}]
The formula \eqref{formpsibetamu} for $\psi_{\beta,\mu}$ shows that $\mu\mapsto \psi_{\beta,\mu}$ is affine and weak* continuous, and both sets of states are weak* compact, so it suffices to show that $\mu\mapsto \psi_{\beta,\mu}$ is surjective and one-to-one. 

To see that $\mu\mapsto \psi_{\beta,\mu}$ is surjective, suppose that $\phi$ is a KMS$_\beta$ state of $\TT(M_L)$. On $C^*(u)=C(\T^d)$, the conditioned state $\phi_P$ is given by a probability measure $\mu$; we choose $\mu$ such that \[
\phi_P(u_m)=\int_{\T^d} z^{m}\,d\mu(z)\ \text{ for $m\in \Z^d$},
\]
and aim to prove that $\phi=\psi_{\beta,\mu}$. Since both states are KMS$_\beta$ states,  formula \eqref{charKMS} shows that it suffices to check that $\phi(u_m)=\psi_{\beta,\mu}(u_m)$. Since $\Z^d$ is abelian, the reconstruction formula \eqref{reconform} implies that 
\begin{align*}
\phi(u_m)&=\lim_{n\to\infty}(1-|\det A|e^{-\beta})\sum_{j=0}^n\sum_{g\in \Z^d/B^j\Z^d} e^{-j\beta}\phi_P(v^{*j}u_mv^j)\\
&=\lim_{n\to\infty}(1-|\det A|e^{-\beta})\sum_{j=0}^n|\det A|^je^{-j\beta}\phi_P(v^{*j}u_mv^j)\\
&=\lim_{n\to\infty}(1-|\det A|e^{-\beta})\sum_{\{j\leq n\;:\;m\in B^j\Z^d\}}|\det A|^je^{-j\beta}\phi_P(u_{B^{-j}m})\\
&=(1-|\det A|e^{-\beta})\sum_{\{j\;:\;m\in B^j\Z^d\}}|\det A|^je^{-j\beta}\int_{\T^d}z^{B^{-j}m}\,d\mu(z),
\end{align*}
which by \eqref{formpsibetamu} is precisely $\psi_{\beta,\mu}(u_m)$. We have now proved surjectivity.

To see that our map is one-to-one, suppose that $\mu$ and $\nu$ are probability measures on $\T^d$ and $\psi_{\beta,\mu}=\psi_{\beta,\nu}$. Write $M_\mu(n)$ for the $n$th moment $\int_{\T^d}z^{n}\,d\mu(z)$ of $\mu$, and fix $m\in \Z^d$. Two applications of \eqref{formpsibetamu} show that
\begin{equation}\label{equatemoments}
\sum_{\{j\;:\;m\in B^j\Z^d\}}|\det A|^je^{-j\beta}M_\mu(B^{-j}m)=\sum_{\{j\;:\;m\in B^j\Z^d\}}|\det A|^je^{-j\beta}M_\nu(B^{-j}m).
\end{equation}
The left-hand side of \eqref{equatemoments} can be rewritten as
\begin{align*}
M_\mu(m&)+\sum_{\{j\;:\;j>0,\;m\in B^j\Z^d\}}|\det A|^je^{-j\beta}M_\mu(B^{-j}m)\\
&=M_\mu(m)+|\det A|e^{-\beta}\sum_{\{j\;:\;j>0,\;m\in B^j\Z^d\}}|\det A|^{j-1}e^{-(j-1)\beta}M_\mu(B^{-(j-1)}B^{-1}m)\\
&=M_\mu(m)+|\det A|e^{-\beta}\sum_{\{j'\;:\;B^{-1}m\in B^{j'}\Z^d\}}|\det A|^{j'}e^{-j'\beta}M_\mu(B^{-j'}B^{-1}m),
\end{align*}
which by \eqref{formpsibetamu} is 
\[
\begin{cases}
M_\mu(m)&\text{if $m$ is not in $B\Z^d$}\\
M_\mu(m)+|\det A|e^{-\beta}\psi_{\beta,\mu}(u_{B^{-1}m})&\text{if $m\in B\Z^d$.}
\end{cases}
\]
If $m$ is not in $B\Z^d$, then \eqref{equatemoments} says precisely that $M_\mu(m)=M_\nu(m)$; if $m\in B\Z^d$, then, since $\psi_{\beta,\mu}(u_{B^{-1}m})=\psi_{\beta,\nu}(u_{B^{-1}m})$, subtracting $|\det A|e^{-\beta}\psi_{\beta,\mu}(u_{B^{-1}m})$ from both sides of \eqref{equatemoments} shows that $M_\mu(m)=M_\nu(m)$. Thus $\mu$ and $\nu$ have the same moments, and are therefore equal.  
\end{proof}

\subsection{Limits of KMS states}\label{seclimits} Proposition~\ref{parametrisation} describes all the KMS$_\beta$ states for $\beta>\beta_c:=\log|\det A|$, and Theorem~\ref{betalogN} says there is exactly one KMS$_{\beta_c}$ state when $A$ is a dilation matrix. So it is natural to ask what we can say about the KMS$_{\beta_c}$ states when $A$ is not a dilation matrix. General results from \cite{bra-rob} suggest that we might be able to find other KMS$_{\beta_c}$ states by taking limits of KMS$_\beta$ states as $\beta\to\beta_c$ from above.

\begin{prop}\label{limitKMS}
Let $\mu\in P(\T^d)$. Then there is a decreasing sequence $\beta_n\to \beta_c$ such that $\{\psi_{\beta_n,\mu}\}$ converges weak* to a state $\psi_\mu$, and then $\psi_\mu$ is a KMS$_{\beta_c}$ state of $(T(M_L),\sigma)$.
\end{prop}

\begin{proof}
Choose any decreasing sequence converging to $\beta_c$, and the weak* compactness of the state space implies that there is a subsequence $\{\beta_n\}$  such that $\{\psi_{\beta_n,\mu}\}$ converges in the weak* topology. Now \cite[Proposition~5.3.23]{bra-rob} implies that the limit $\psi_\mu$ is a KMS$_{\beta_c}$ state, at least when $\beta_c>0$. When $\beta_c=0$, \cite[Proposition~5.3.23]{bra-rob} only asserts that $\psi_\mu$ is a trace (because that is what being a KMS$_0$ state means in \cite{bra-rob}). However, KMS$_\beta$ states for $\beta>0$ are $\sigma$-invariant, and hence so is the limit. Thus $\psi_\mu$ is a KMS$_0$ state in the sense we are using. 
\end{proof}

We now assume that $A$ is not a dilation matrix, so that $\bigcap_{j=0}^\infty B^j\Z^d$ could be bigger than $\{0\}$. Suppose $\beta>\beta_c$ and write $r=e^{-(\beta-\beta_c)}$. As in the last proof, we write $M_\mu(m)$ for the $m$th moment $\int_{\T^d}z^m\,d\mu(z)$. Rearranging \eqref{formpsibetamu} shows that $\psi_{\beta,\mu}(u_mv^kv^{*l}u_n^*)$ vanishes unless $k=l$ and $m-n\in B^k\Z^d$, and then equals
\begin{equation}\label{formrewritten}
\sum_{\{j\geq 0\,:\,m-n\in B^{j+k}\Z^d\}}e^{-k\beta}(1-r)r^jM_\mu(B^{-(j+k)}(m-n)).
\end{equation}
So we want to compute the limit of \eqref{formrewritten} as $\beta\to \beta_c$, in which case $r\to 1-$. If $m-n$ does not belong to $\bigcap_{j=0}^\infty B^{j+k}\Z^d=\bigcap_{j=0}^\infty B^{j}\Z^d$, then the sum in \eqref{formrewritten} is finite, and since $(1-r)r^j\to 0$ as $r\to 1$ for each fixed $j$, \eqref{formrewritten} converges to $0$ as $r\to 1$. So it remains for us to compute the limit of \eqref{formrewritten} when $m-n\in \bigcap_{j=0}^\infty B^{j}\Z^d$. Unfortunately, this seems to be a fairly delicate matter (see Remark~\ref{cantuseDCT} below), and the best we can do is illustrate the issues with some examples.

\begin{itemize}
\item[(a)] If $\mu$ is normalised Haar measure on $\T^d$, then $M_\mu(0)=1$ and $M_\mu(m)=0$ for all other $m$. The series in \eqref{formrewritten} is identically zero unless $m=n$, and then is geometric; summing it shows that 
\[
\psi_{\beta,\mu}(u_mv^kv^{*l}u_n^*)=
\begin{cases}
0&\text{unless $k=l$ and $m=n$}\\
e^{-k\beta}&\text{if $k=l$ and $m=n$.}
\end{cases}
\]
Letting $\beta\to \beta_c$ gives the state described in Theorem~\ref{betalogN}.

\item[(b)] If $\mu$ has the property that $M_\mu(m)=1$ for every $m\in\bigcap_{j=0}^\infty B^{j}\Z^d$, then the series in \eqref{formrewritten} is geometric whenever $m-n\in \bigcap_{j=0}^\infty B^{j}\Z^d$. Summing and letting $\beta\to\beta_c$ shows that the limit $\psi_\mu$ satisfies
\[
\psi_{\mu}(u_mv^kv^{*l}u_n^*)=
\begin{cases}
0&\text{unless $k=l$ and $m-n\in\bigcap_{j=0}^\infty B^{j}\Z^d$}\\
e^{-k\beta_c}&\text{if $k=l$ and $m-n\in\bigcap_{j=0}^\infty B^{j}\Z^d$.}
\end{cases}
\]
\item[(c)] The previous item (b) applies in particular to the point mass $\delta_1$ at the identity $1=(1,1,\cdots,1)$ of $\T^d$. This shows that the KMS$_{\beta_c}$ state in Theorem~\ref{betalogN} is unique if and only if  $\bigcap_{j=0}^\infty B^{j}\Z^d=\{0\}$.

\item[(d)] Consider the matrix $A=\big(\begin{smallmatrix}2&0\\0&1\end{smallmatrix}\big)$, for which $\bigcap_{j=0}^\infty B^{j}\Z^2=\{0\}\times\Z$. Then item (b) applies to any measure of the form $\nu\times\delta_1$. Thus when $A$ is not a dilation matrix, we expect there to be many KMS$_{\beta_c}$ states besides the one in Theorem~\ref{betalogN}.

\item[(e)] We wonder whether every KMS$_{\beta_c}$ state is a limit of KMS$_\beta$ states. It is trivially the case in our examples when $|\det A|>1$, and in these examples it also works for $\beta_c=0$. 

\item[(f)] When $\beta_c=0$, we have to be careful to distinguish between traces (the KMS$_0$ states in \cite{bra-rob}) and the invariant traces (the KMS$_0$ states in \cite{ped}). Certainly any limit of KMS$_\beta$ states will be invariant, so the answer to the previous question is trivially false with the definition in \cite{bra-rob} if the algebra has traces which are not invariant. We give an example where this happens in Remark~\ref{tracesonHgp}.
\end{itemize}

\begin{remark}\label{cantuseDCT} The obvious way to try to compute the limit of \eqref{formrewritten} as $r\to 1-$ is to evaluate it term-by-term. This amounts to pulling $\lim_{r\to 1-}$ through the infinite sum, and therefore requires the dominated convergence theorem. Write $m_j:=M_\mu(B^{-(j+k)}(m-n))$. To apply the dominated convergence theorem, we need a convergent series $\sum_j a_j$ such that $0\leq (1-r)r^j|m_j|\leq a_j$ (and we need to consider a sequence $\{r_n\}$). We know $|m_j|\leq 1$. Calculus shows that $\max\{(1-t)t^j:t\in[0,1]\}$ occurs at $j/(j+1)$. So the best general estimate seems to be 
\[(1-r)r^j|m_j|\leq \Big(1-\frac{j}{j+1}\Big)\Big(\frac{j}{j+1}\Big)^j=\frac{j^j}{(j+1)^{j+1}}.\]
Taking $a_j$ to be the right-hand side and $b_j:=1/(j+1)$, we have
\[\frac{b_j}{a_j}=\Big(\frac{j+1}{j}\Big)^j=\Big(1+\frac{1}{j}\Big)^j\to e\ \text{ as $j\to \infty$,}\]
and the limit form of the comparison test implies that $\sum a_j$ diverges. 

So pulling the limit through the sum seems to be a nontrivial matter. Of course, it is really just as well we can't do this, since we know that $\sum_{j=0}^\infty (1-r)r^j =1 \to 1$ as $r\to 1-$, whereas the term-by-term calculation would give $0$.
\end{remark}

\section{KMS$_\infty$ and ground states}\label{sec-ground}

\begin{prop}\label{betainfty}
Suppose that $A\in M_d(\Z)$ has nonzero determinant. Then for every probability measure $\mu$ on $\T^d$, there is a KMS$_\infty$ state $\psi_{\infty,\mu}$ on $(\TT(M_L),\sigma)$ such that
\[
\psi_{\infty,\mu}(u_mv^kv^{*l}u_n^*)=
\begin{cases}\int_{\T^d}z^{m-n}\,d\mu(z)&\text{if $k=l=0$}\\
0&\text{otherwise.}
\end{cases}
\]
Every ground state of $(\TT(M_L),\sigma)$ has the form $\psi_{\infty,\mu}$, and is in particular a KMS$_\infty$ state.
\end{prop}

The proof of \cite[Lemma~8.4]{lr} gives the following characterisation of ground states.

\begin{lemma}\label{charground}
A state $\phi$ of $\TT(M_L)$ is a ground state for $\sigma$ if and only if
\[
\phi_{\infty,\mu}(u_mv^kv^{*l}u_n^*)=
\begin{cases}\phi(u_{m-n})&\text{if $k=l=0$}\\
0&\text{otherwise.}
\end{cases}
\]
\end{lemma}

\begin{proof}[Proof of Proposition~\ref{betainfty}]
Choose a sequence $\{\beta_i\}$ such that $\beta_i\to \infty$; by passing to a subsequence, we may suppose that $\psi_{\beta_i,\mu}$ converges in the weak* topology to a state $\psi_{\infty,\mu}$, which is by definition a KMS$_\infty$ state. Next we verify the formula for $\psi_{\infty,\mu}$. As $\beta\to \infty$, each summand in the right-hand side of \eqref{formpsibetamu} with $j>0$ goes to zero. Thus as $i\to \infty$, we have
\[
\psi_{\beta_i,\mu}(u_mv^kv^{*l}u_n^*)\to 
\begin{cases}\int_{\T^d}z^{m-n}\,d\mu(z)&\text{if $k=l=0$}\\
0&\text{otherwise,}
\end{cases}
\]
and hence $\phi_{\infty,\mu}$ has the required form. (If $A$ is not a dilation matrix, so that $\bigcap_{k=1}^\infty B^k\Z^d$ could contain nonzero elements, then the sum on the right-hand side of \eqref{formpsibetamu} could have infinitely many nonzero terms, and calculating the limit as $i\to \infty$ would require some analysis of the sort discussed in Remark~\ref{cantuseDCT}.) 

If $\phi$ is a ground state, then the restriction of $\phi$ to the range of $i_{C(\T^d)}$ is given by a probability measure $\mu$, and then Lemma~\ref{charground} implies that $\phi=\psi_{\infty,\mu}$.
\end{proof}

\section{The Toeplitz algebra of the Baumslag-Solitar semigroup}\label{BSsemigp}

We fix an integer $N$ with $N>1$, and consider the additive group $\Z[N^{-1}]$ of rational numbers of the form $mN^{-l}$ for $m,l\in \Z$. The \emph{Baumslag-Solitar group} is the semidirect product $\Z[N^{-1}]\rtimes \Z$ with
\[
(r,k)(s,l)=(r+N^ks, k+l).
\]

The semigroup semidirect product $\N\rtimes_N\N$ is a subsemigroup of $\Z[N^{-1}]\rtimes \Z$, and the pair $(\Z[N^{-1}]\rtimes \Z,\N\rtimes_N \N)$ is closely related to the pair $(\qxqx,\nxnx)$ studied in \cite{lr}. Indeed, the map $(r,k)\mapsto (r,N^k)$ of $\Z[N^{-1}]\rtimes\N$ into $\qxqx$ carries $\N\rtimes_N\N$ into $\nxnx$. The pair $(\Z[N^{-1}]\rtimes \Z,\N\rtimes_N \N)$ is also quasi-lattice ordered in the sense of Nica \cite{n}. (One way to see this is via the embedding of $\Z[N^{-1}]\rtimes\N$ in $\qxqx$: we just need to check that if $(r,k)\in\Z[N^{-1}]\rtimes\N$ and $(r,N^k)$ has an upper bound in $\N\rtimes N^\N$, then the least upper bound constructed in \cite[Proposition~2.2]{lr} lies in $\N\rtimes N^\N$.) So $(\Z[N^{-1}]\rtimes \Z,\N\rtimes_N \N)$ also has a Toeplitz algebra $\TT(\N\rtimes_N\N)$ which is universal for Nica covariant representations of  $\N\rtimes_N\N$ \cite{n,lr1}. The Toeplitz algebra $\TT(C(\T),\alpha_N,L,\N)$ is a quotient of $\TT(\N\rtimes_N\N)$ analogous to the additive boundary quotient of $\TT(\nxnx)$ studied in \cite{BaHLR}. We now discuss the KMS states on $\TT(\N\rtimes_N\N)$, following the analysis of \cite[\S4]{BaHLR}.

The Toeplitz algebra $\TT(\N\rtimes_N \N)$ is generated by two isometries $s=T_{(1,0)}$ and $v=T_{(0,1)}$, and an argument like that of \cite[\S4]{lr} shows that $(\TT(\N\rtimes_N\N),s,v)$ is universal among $C^*$-algebras generated by a pair of isometries $S$ and $V$ satisfying
\begin{itemize}
\item[(T1)]\  $VS =  S^NV$,
\smallskip
\item[(T4)]\ $S^*V = S^{N-1}V S^*$, and
\smallskip
\item[(T5)]\ $V^* S^k V = 0$ for $1 \leq k < N$.
\end{itemize}
We define $\TT_{\add}(\N\rtimes_N \N)$ to be the quotient of $\TT(\N\rtimes_N \N)$ by the extra relation $ss^*=1$, $\TT_{\mult}(\N\rtimes_N \N)$ to be the quotient by the relation $1=\sum_{k=0}^{N-1}s^kvv^*s^{*k}$, and $\TT_{\add,\;\mult}(\N\rtimes_N \N)$ to be the quotient in which both extra relations hold, and which is therefore the analogue of Cuntz's $\qn$. Thus we have the following commutative diagram of quotient maps:
\begin{equation}\label{commdiagBS}
\xymatrix{
&\TT(\N\rtimes_N \N)\ar[dl]_{q_{\add}}\ar[dr]^{q_{\mult}}&\\
\TT_{\add}(\N\rtimes_N \N)\ar[dr]&&\TT_{\mult}(\N\rtimes_N \N)\ar[dl]\\
&\TT_{\add,\;\mult}(\N\rtimes_N \N).&
}
\end{equation}

In $\TT_{\add}(\N\rtimes_N \N)$ the generator $s$ becomes unitary, and (T4) is redundant. The unitary $s$ generates a unitary representation $u:\Z\to U(\TT_{\add}(\N\rtimes_N \N))$, and the relations (T1) and (T5) (taken together) are equivalent to (E1) and (E2) (taken together). Thus Proposition~\ref{presadd} implies that $\TT_{\add}(\N\rtimes_N \N)$ is our Toeplitz algebra $\TT(M_L)$. Proposition~\ref{CPquotient} implies that $\TT_{\add,\;\mult}(\N\rtimes_N \N)$ is the quotient $\OO(M_L)$ of $\TT(M_L)$. (When $N=2$, $\OO(M_L)$ has been studied by Larsen and Li under the name $\mathcal{Q}_2$; see \cite[\S3]{LL}.)

Since the presentation of $\TT(\N\rtimes_N\N)$ is not affected by multiplying $v$ by $z\in \T$, we can deduce from the presentation that there is an action $\gamma: \T\to \Aut\TT(\N\rtimes_N\N)$ such that $\gamma_z(s)=s$ and $\gamma_z(v)=zv$. Inflating this action to $\R$ gives a  dynamics $\sigma:\R\to\Aut\TT(\N\rtimes_N\N)$ such that $\sigma_t(s)=s$ and $\sigma_t(v)=e^{it}v$. This action leaves the kernels of the quotient maps in the diagram \eqref{commdiagBS} invariant, and hence induces actions (still denoted by $\sigma$) on all three quotients. On $\TT_{\add}(\N\rtimes_N \N)$ and $\TT_{\add,\;\mult}(\N\rtimes_N \N)$ we recover the actions on $\TT(M_L)$ and $\OO(M_L)$ that we have been studying, in the case where $A$ is the $1\times 1$ matrix $(N)$ and $\sigma_A$ is the covering map $z\mapsto z^N$ of $\T$. So our results tell us about the KMS states of $(\TT_{\add}(\N\rtimes_N \N),\sigma)$ and $(\TT_{\add,\;\mult}(\N\rtimes_N \N),\sigma)$.

Just as in \cite[Lemma~10.4]{lr}, every KMS$_\beta$ state of $(\TT(\N\rtimes_N\N),\sigma)$ vanishes on the ideal generated by $1-ss^*$, and hence comes from a KMS$_\beta$ state of $\TT(M_L)$. So we know all the KMS states of $(\TT(\N\rtimes_N\N),\sigma)$. For ground states, though, there is a difference. As in \cite[Lemma~8.4]{lr} (or Lemma~\ref{charground} above), a ground state of $(\TT(\N\rtimes_N\N),\sigma)$ is determined by its values on $C^*(s)$, and we claim that the map $\phi\mapsto \phi|_{C^*(s)}$ is an affine homeomorphism of the set of ground states onto the state space of $C^*(s)\cong \TT(\N)$. Indeed, we can deduce this from \cite[Theorem~7.1(4)]{lr}, since Theorem~3.7 of \cite{lr1} implies that $\TT(\N\rtimes_N\N)$ embeds as the subalgebra $C^*(s, v_N)$ of $\TT(\nxnx)$, and the homeomorphism $\phi\mapsto \phi|_{C^*(s)}$ factors through $C^*(s,v_N)$.

We can sum up these results by saying that the system $(\TT(\N\rtimes_N\N),\sigma)$ has a phase transition at inverse temperature $\beta=\log N$, and a further phase transition (in the sense of Connes and Marcolli) at $\beta=\infty$. We believe that this is the simplest known system which exhibits both these phenomena. As for the system in \cite{lr}, the circular symmetry at $\beta=\log N$ which disappears for $\beta>\log N$ is not apparently realised by an action of $\T$ on $\TT(\N\rtimes_N\N)$. In \cite{lr}, though, this circular symmetry persists for $\beta\in [1,2]$, as a result of the more complicated convergence issues for the series representations of the normalising factors. 

\begin{remark}
Since we can view $\N\rtimes_N\N$  as a subsemigroup of $\nxnx$, it might be more natural to use the dynamics satisfying $\sigma_t(v)=N^{it}v$. If we do this, then the phase transition will occur at $\beta=1$.
\end{remark}

\section{Integer matrices with determinant $\pm 1$}\label{sec-gpalgs}

When $A\in M_d(\Z)$ has $|\det A|=1$, the inverse $A^{-1}$ has integer entries (as the cofactor formula shows), the map $\sigma_A$ is a homeomorphism, and $\alpha_A$ is an automorphism. The inverse $\alpha_A^{-1}$ is then a transfer operator for $\alpha_A$, so we have an Exel system $(C(\T^d),\alpha_A,\alpha_A^{-1})$, and this system has a Toeplitz algebra and an Exel crossed product. One would guess that these $C^*$-algebras must be related to the ordinary crossed product, and they are, but we have not seen this explicitly pointed out before. 

\begin{prop}
Suppose that $\alpha$ is an automorphism of a unital $C^*$-algebra $C$. Then

\smallskip
\textnormal{(a)} Exel's Toeplitz algebra $\TT(C,\alpha,\alpha^{-1})$ is the universal $C^*$-algebra generated by an isometry $v$ and a unital representation $i_C$ of $C$ satisfying $vi_C(c)=i_C(\alpha(c))v$, and

\smallskip
\textnormal{(b)} the Exel crossed product $C\rtimes_{\alpha,\alpha^{-1}}\N$ is the universal $C^*$-algebra generated by a unitary $u$ and a unital representation $j_C$ of $C$ satisfying $j_C(\alpha(c))=uj_C(c)u^*$.
\end{prop}

\begin{proof}
We know from \cite[\S3]{br} that $\TT(C,\alpha,\alpha^{-1})$ is universal for Toeplitz-covariant representations $(\rho, V)$ satisfying two relations called (TC1) and (TC2) (see page~\pageref{TCrels}). As we observed earlier, plugging the identity $1$ of $C$ into (TC2) shows that $V$ is an isometry. For our system  $(C,\alpha,\alpha^{-1})$, (TC1) implies (TC2):
\[
V^*\rho(c)V=V^*\rho(\alpha(\alpha^{-1}(c)))V=V^*V\rho(\alpha^{-1}(c))=\rho(\alpha^{-1}(c)),
\]
and (a) follows. 

To establish (b), notice first that $\phi(c)\in\LL(M_{\alpha^{-1}})$ is the rank-one operator $\Theta_{c,1}$. The Cuntz-Pimsner algebra is generated by a universal Cuntz-Pimsner covariant representation $(j_{M_{\alpha^{-1}}},j_C)$, and then the isometry $v$ in part (a) is $v=j_{M_{\alpha^{-1}}}(1)$. Cuntz-Pimsner covariance says that
\[
j_C(c)=(j_{M_{\alpha^{-1}}},j_C)^{(1)}(\phi(c))=(j_{M_{\alpha^{-1}}},j_C)^{(1)}(\Theta_{c,1})=j_{M_{\alpha^{-1}}}(c)j_{M_{\alpha^{-1}}}(1)^*;
\]
since $c=c\cdot 1$, we have $j_{M_{\alpha^{-1}}}(c)=j_C(c)j_{M_{\alpha^{-1}}}(1),$, and Cuntz-Pimsner covariance is equivalent to $j_C(c)=j_C(c)vv^*$. This is equivalent to $vv^*=1$, so $v$ is unitary, and now $vj_C(c)=j_C(\alpha(c))v$ is equivalent to $j_C(\alpha(c))=vj_C(c)v^*$.
\end{proof}

These universal properties immediately imply that our algebras are familiar objects:

\begin{cor}\label{ordcp}
Suppose that $\alpha$ is an automorphism of a unital $C^*$-algebra $C$. Then the Exel crossed product $C\rtimes_{\alpha,\alpha^{-1}}\N$ is the usual crossed product $C\rtimes_\alpha\Z$, and the Toeplitz algebra $\TT(C,\alpha,\alpha^{-1})$ is the crossed product $C\times_\alpha\N$ introduced and studied by Murphy \cite{m}. In both cases, the gauge action of $\T$ is the dual action of $\T=\widehat\Z$.
\end{cor}

\begin{remark}
Although $M_{\alpha^{-1}}$ is not the bimodule $E$ considered by Pimsner in \cite[Example~(3), page~193]{p}, the two are isomorphic; indeed, $a\mapsto \alpha(a)$ is a Hilbert-bimodule isomorphism of $E$ onto $M_{\alpha^{-1}}$. So the identity $C\rtimes_{\alpha,\alpha^{-1}}\N=C\rtimes_\alpha\Z$ also follows from the assertion in \cite[Example~(3)]{p}.
\end{remark}

We now return to the case of an integer matrix $A$ with $|\det A|=1$, where Cor\-ollary~\ref{ordcp} identifies the Toeplitz algebra $\TT(C(\T^d),\alpha_A,\alpha_A^{-1})$ as a Murphy crossed product, and the Exel crossed product $C(\T^d)\rtimes_{\alpha_A,\alpha_A^{-1}}\N$ with the ordinary crossed product $C(\T^d)\rtimes_{\alpha_A}\Z$. We will be working primarily with the crossed product, so it is worth observing that the generator $\bar v$ is now unitary, and hence we can simplify our presentation: we view $C(\T^d)\rtimes_{\alpha}\Z$ as being generated by a unitary representation $\bar u$ of $\Z^d$ and a unitary $\bar v$ satisfying $\bar v\bar u_m \bar v^*=\bar u_{Bm}$, and then
\[
C(\T^d)\rtimes_{\alpha_A}\Z=\clsp\{\bar u_m\bar v^k: m\in \Z^d, k\in\Z\}.
\]
We can if we wish make the further identification of the crossed product $C(\T^d)\rtimes_{\alpha_A}\Z$ with the group algebra $C^*(\Z^d\rtimes_B\Z)$ of the semidirect product (using Proposition 3.11 of \cite{tfb^2}, for example).

As before, lifting the dual actions of $\T$ gives actions $\sigma$ of $\R$ on $\TT(C(\T^d),\alpha_A,\alpha_A^{-1})$ and $C(\T^d)\rtimes_{\alpha_A}\Z$ such that $\sigma_t$ fixes the copies of $C(\T^d)$ and multiplies the additive generators by $e^{it}$. Proposition~\ref{existforbeta>1} describes the KMS$_\beta$ states of $(\TT(C(\T^d),\alpha_A,\alpha_A^{-1}),\sigma)$ for $\beta>\log|\det A|=0$. Since $|\det B|=1$, $B$ is invertible over the integers,  $B^j\Z^d=\Z^d$ for all $j$, and the series in \eqref{formpsibetamu} is infinite for every pair $m$, $n$. Thus for each $\mu\in P(\T^d)$ there is a KMS$_\beta$ state $\psi_{\beta,\mu}$ on $\TT(C(\T^d),\alpha_A,\alpha_A^{-1})$ such that
\begin{equation}\label{psibeta}
\psi_{\beta,\mu}(u_mv^kv^{*l}u_n^*)=
\begin{cases}
0&\text{unless $k=l$}\\
\sum_{j=k}^\infty (1-e^{-\beta})e^{-j\beta}M_\mu(B^{-j}(m-n))&\text{if $k=l$.}
\end{cases}
\end{equation}
Indeed, the proof of Proposition~\ref{existforbeta>1} simplifies substantially in this case: with $U=\bigoplus_{j=0}^\infty M\circ B^{-j}$ acting on $\bigoplus_{j=0}^\infty L^2(\T^d,d\mu)$, $V$ the unilateral shift on the same direct sum, and $e_j$ the constant function $1$ in the $j$th summand and $0$ elsewhere, we have
\[
\psi_{\beta,\mu}(T)=\sum_{j=0}^\infty(1-e^{-\beta})e^{-j\beta}(\pi_{U,V}(T)e_j\,|\,e_j).
\]
Proposition~\ref{existforbeta>1} also shows that all the KMS$_0$ states (that is, the invariant traces) on $\TT(C(\T^d),\alpha_A,\alpha_A^{-1})$ factor through traces of $C^*(\T^d)\rtimes_{\alpha_A}\Z$. Since the uniqueness assertion in Theorem~\ref{betalogN} does not apply, we might expect to find more than one.

\begin{prop}\label{invtr}
Suppose that $A\in M_d(\Z)$ has $|\det A|=1$. If $\mu\in P(\T^d)$ satisfies $\sigma_A^*\mu=\mu$, then there is a $\sigma$-invariant trace $\psi_\mu$ on $C(\T^d)\rtimes_{\alpha_A}\Z$ such that
\begin{equation}\label{proppsimu}
\psi_\mu(\bar u_m\bar v^k)=
\begin{cases}
0&\text{unless $k=0$}\\ M_m(\mu)&\text{if $k=0$,}
\end{cases}
\end{equation}
 and every $\sigma$-invariant trace on $C(\T^d)\rtimes_{\alpha_A}\Z$ has this form.
\end{prop}

Since the action $\sigma$ of $\R$ is inflated from the dual action $\widehat\alpha_A$ of $\T$, a state is invariant for $\sigma$ if and only if it is invariant for $\widehat\alpha_A$. So the following standard lemma is useful.

\begin{lemma}\label{invstate}
Suppose that $\gamma:\T\to \Aut D$ is a strongly continuous action on a unital $C^*$-algebra $D$ and $E^\gamma:d\mapsto \int_{\T} \gamma_z(d)\,dz$ is the expectation onto the fixed-point algebra $D^\gamma$. Then a state $\phi$ of $D$ is $\gamma$-invariant if and only if there is a state $\tau$ of $D^\gamma$ such that $\phi=\tau\circ E^\gamma$. 
\end{lemma}

\begin{proof}
Suppose $\phi=\tau\circ E^\gamma$. Then the invariance of Haar measure implies that $E^\gamma\circ \gamma_z=E^\gamma$, and hence $\phi\circ \gamma_z=\tau\circ E^\gamma\circ\gamma_z=\tau\circ E^\gamma=\phi$, so $\phi$ is invariant. Conversely, if $\phi$ is invariant, then 
\[
\phi(d)=\int_{\T}\phi(d)\,dz=\int_{\T}\phi(\gamma_z(d))\,dz=\phi\Big(\int_{\T}\gamma_z(d)\,dz\Big)=\phi\circ E^\gamma(d),
\]
so $\phi=\phi|_{D^\gamma}\circ E^\gamma$; since $1\in D^\gamma$, $\phi|_{D^\gamma}$ is a state.
\end{proof}

\begin{proof}[Proof of Proposition~\ref{invtr}]
With $\gamma=\widehat\alpha_A$, the expectation $E^{\gamma}$ is given by
\[
E^\gamma(\bar u_m\bar v^k)=\begin{cases}
0&\text{unless $k=0$}\\ \bar u_m&\text{if $k=0$.}
\end{cases}
\]
It follows easily that if $\theta_\mu$ is the state on $C(\T^d)=\clsp\{\bar u_m:m\in\Z^d\}$ given by integration against $\mu\in P(\T^d)$, then $\psi_\mu:=\theta_\mu\circ E^\gamma$ satisfies \eqref{proppsimu}. Lemma~\ref{invstate} implies that every invariant state of $C(\T^d)\rtimes_{\alpha_A}\Z$ has this form. So we need to show that $\psi_\mu$ is a trace if and only if $\mu$ is invariant under $\sigma_A^*$. 

We compute
\[
\psi_\mu((\bar u_m\bar v^k)(\bar u_n\bar v^l))
=\begin{cases}
\psi_\mu(\bar u_{m+B^kn})&\text{if $k+l=0$}\\
0&\text{otherwise,}
\end{cases}
\]
and similarly the other way round.
Thus $\psi_\mu$ is a trace if and only if 
\[
\psi_\mu(\bar u_{m+B^{-l}n})=\psi_\mu(\bar u_{n+B^{l}m})=\psi_\mu(\bar u_{B^l(B^{-l}n+m)})\ \text{ for all $l\in\Z$, $m,n\in\Z^d$;}
\]
or, equivalently, if and only if $\psi_\mu(\bar u_m)=\psi_\mu(\bar u_{Bm})$ for all $m\in \Z^d$. But these are just the moments of $\mu$, and a calculation shows that $M_{Bm}(\mu)=M_m(\sigma_A^*\mu)$. So we deduce that $\psi_\mu$ is a trace if and only if $\mu$ is invariant, as required.
\end{proof}

\begin{remark}
When $\mu$ is invariant under $\sigma_A$, the moments $M_\mu(B^{-j}(m-n))$ appearing in \eqref{psibeta}  are all equal to $M_\mu(m-n)$. Thus the series on the right-hand side of \eqref{psibeta} is geometric, and summing it shows that
\begin{equation*}
\psi_{\beta,\mu}(u_mv^kv^{*l}u_n^*)=
\begin{cases}
0&\text{unless $k=l$}\\
e^{-k\beta}M_\mu(m-n)&\text{if $k=l$,}
\end{cases}
\end{equation*}
which converges to $\psi_\mu(Q(u_mv^kv^{*l}u_n^*))=\psi_\mu(\bar u_{m+B^{k-l}n}\bar v^{k-l})$ as $\beta\to 0$. In view of Proposition~\ref{limitKMS}, this gives an alternative proof that $\psi_\mu$ is an invariant trace.
\end{remark}

\begin{cor}
When $|\det A|=1$, $C(\T^d)\rtimes_{\alpha_A}\Z$ has many invariant traces.
\end{cor}

\begin{proof}
The homeomorphism $\sigma_A$ has many finite orbits --- indeed, the periodic points are dense in $\T^d$. (If the coordinates of $r\in \Q^d$ have common denominator $N$, then the denominators of all the coordinates in all the $A^nr$ divide $N$ too, and the pigeon-hole principle implies that $e^{2\pi ir}$ is periodic.) But if $z\in \T^d$ has $\sigma_A^n(z)=z$, then $n^{-1}\sum_{j=0}^{n-1}\delta_{\sigma_A^j(z)}$ is an invariant measure.
\end{proof}

\begin{example}\label{tracesHgp}
Consider $A=\big(\begin{smallmatrix}1&0\\1&1\end{smallmatrix}\big)$, for which we have $\sigma_A(w,z)=(w,wz)$. Let $\lambda$ denote Haar measure on $\T$ and let $\nu$ be a probability measure on $\T$. Then the product measure $\mu=\nu\times\lambda$ is invariant for $\sigma_A$:
\begin{align*}
\int_{\T^2}f\circ\sigma_A&(w,z)\,d\mu(w,z)=\int_{\T}\!\int_{\T}f(w,wz)\,d\lambda(z)\,d\nu(w)\\
&=\int_{\T}\!\int_{\T}f(w,z)\,d\lambda(z)\,d\nu(w)=\int_{\T^2}f(w,z)\,d\mu(w,z).
\end{align*}
The moments of $\mu$ are given by
\[
M_\mu(m)=\int_{\T}\!\int_{\T} w^{m_1}z^{m_2}\,d\lambda(z)\,d\nu(w)=
\begin{cases}
0&\text{unless $m_2=0$}\\
M_\nu(m_1)&\text{if $m_2=0$.}
\end{cases}
\]
Thus Proposition~\ref{invtr} gives $\sigma$-invariant traces $\{\phi_\nu:\nu \in P(\T)\}$ on $C^*(\T^d)\rtimes_{\alpha_A}\Z$ such that
\[
\phi_\nu(\bar u_m\bar v^k)=
\begin{cases}
0&\text{unless $k=0$ and $m_2=0$}\\
M_\nu(m_1)&\text{if $k=0$ and $m_2=0$.}
\end{cases}
\]
\end{example}

\begin{example}
When $A=\big(\begin{smallmatrix}1&0\\1&1\end{smallmatrix}\big)$ and $\mu=\lambda\times\nu$, we claim that the state $\psi_{\beta,\mu}$ described in \eqref{psibeta} converges as $\beta\to 0+$ to the state of Theorem~\ref{betalogN}, which for $\beta=\log|\det A|=0$ vanishes on $u_mv^k$ unless $(m,k)=(0,0)$ and satisfies $\psi(u_0)=1$. To see this, note that $B=A^t=\big(\begin{smallmatrix}1&1\\0&1\end{smallmatrix}\big)$, so $B^{-j}=\big(\begin{smallmatrix}1&-j\\0&1\end{smallmatrix}\big)$, and
\[
M_{\lambda\times\nu}(B^{-j}(m-n))=\int\!\int w^{(m_1-n_1)-j(m_2-n_2)}\,d\lambda(w)\,z^{m_2-n_2}\,d\nu(z),
\]
which vanishes unless $m_1-n_1=j(m_2-n_2)$. For $m=n$, we have $M_\mu(B^{-j}(m-n))=1$ for all $j$, and summing the series shows that $\psi_{\beta,\mu}(u_nv^kv^{*k}u_n^*)=e^{-k\beta}\to 1$ as $\beta\to 0$ for all $k$. If $m\not=n$ and $m-n$ does not have the form $(li,i)$, then $\psi_{\beta,\mu}(u_mv^kv^{*k}u_n^*)=0$ for all $k$. If $m-n=(li,i)$, then $\psi_{\beta,\mu}(u_mv^kv^{*k}u_n^*)$ vanishes for $k>l$, and
\[
\psi_{\beta,\mu}(u_mv^kv^{*k}u_n^*)=(1-e^{-\beta})e^{-l\beta}M_\nu(i)\to 0\ \text{ as $\beta\to 0$}
\]
for $k\leq l$. So whenever $m\not=n$, we have $\psi_{\beta,\mu}(u_mv^kv^{*k}u_n^*)\to 0$ as $\beta\to 0$, and the limit is the state $\phi$ of Theorem~\ref{betalogN}, as claimed. (Well, strictly speaking the limit is $\phi\circ Q$.)
\end{example}

\begin{remark}\label{tracesonHgp}
When $A=\big(\begin{smallmatrix}1&0\\1&1\end{smallmatrix}\big)$ and $B=A^t=\big(\begin{smallmatrix}1&1\\0&1\end{smallmatrix}\big)$, the map 
\begin{equation}\label{defHgp}(m,k)=((m_1,m_2),k)\mapsto \begin{pmatrix}1&k&m_1\\0&1&m_2\\0&0&1
\end{pmatrix}
\end{equation}
is an isomorphism of $\Z^2\rtimes_B\Z$ onto the integer Heisenberg group $H(\Z)$. (The crux is that $B^k=\big(\begin{smallmatrix}1&k\\0&1\end{smallmatrix}\big)$.) For $\theta\in [0,1]$, we view the rotation algebra $\AA_\theta$ as the universal $C^*$-algebra generated by unitaries $U,V$ satisfying $VU=e^{2\pi i\theta}UV$, and then the unitary representation $(m,k)\mapsto e^{2\pi im_1\theta}U^{m_2}V^k$ induces a surjection $q_\theta$ of $C^*(H(\Z))$ onto $\AA_\theta$. (Indeed, the quotients $\AA_\theta$ are the fibres of a $C^*$-bundle over $\T$ which has $C^*(H(\Z))$ as its algebra of continuous sections --- see \cite{e}, \cite[\S1]{ap} or \cite[Example~1.4]{prae}.)

Every rotation algebra $\AA_\theta$ has a trace $\tau_\theta$ which kills $U^{m_2}V^k$ unless $m_2=0=k$, and the composition $\tau_\theta\circ q_\theta$ is the invariant trace described in Example~\ref{tracesHgp} for $\nu$ the point mass at $e^{2\pi i\theta}$. When $\theta$ is irrational, $\tau_\theta$ is the only trace on $\AA_\theta$ (see \cite[Proposition~VI.1.3]{dav}, for example). When $\theta$ is rational, $\AA_\theta$ is a homogeneous $C^*$-algebra with spectrum $\T^2$ (by, for example, \cite[\S2]{dr}), and has other traces which give non-invariant traces of $C^*(H(\Z))$. For example, the matrices $T=\big(\begin{smallmatrix}0&1\\1&0\end{smallmatrix}\big)$ and $S=\big(\begin{smallmatrix}-1&0\\0&1\end{smallmatrix}\big)$ are unitary and satisfy $TS=-ST$, hence give a homomorphism $\pi_{S,T}:\AA_{1/2}\to M_2(\C)$, and composing with the usual normalised trace $2^{-1}\tr$ gives a trace $\tau$ on $\AA_{1/2}$ such that $\tau(U^{m_2}V^k)=2^{-1}\tr(S^{m_2}T^k)$. Since $T^2=1$, we have $\tau(V^2)=\tau(1)=1$, and since $\sigma_t(v^2)=e^{2it}v^2$, $\tau(V^2)=1$ implies that $\tau\circ q_{1/2}$ cannot be $\sigma$-invariant.

To explain where the other invariant traces in Example~\ref{tracesHgp} come from, we examine the structure of $C^*(H(\Z))$ from another point of view. Consider the normal subgroup $N$ of matrices \eqref{defHgp} with $k=m_2=0$, which is the centre of $H(\Z)$, and which has quotient $H(\Z)/N$ isomorphic to $\Z^2$ via $(m,k)\mapsto (k,m_2)$. Applying Theorem~4.1 of \cite{prae2} to $N$ gives a realisation of $C^*(H(\Z))=\C\times_{\id,1}H(\Z)$ as a Busby-Smith twisted crossed product $C^*(N)\rtimes_{\beta,\omega}\Z^2$; identifying $C^*(N)$ with $C(\T)$ and ploughing through the formulas in \cite{prae2} shows that $\beta$ is the identity and the cocycle $\omega:\Z^2\to U(C(\T))=C(\T,\T)$ is given by
\[
\omega((k,m_2),(l,n_2))(z)=z^{kn_2}.
\]
Averaging over the dual action of $\T^2$ gives an expectation $E^{\widehat\beta}$ whose range is the fixed-point algebra $C(\T)\subset C(\T)\rtimes_{\id,\omega}\Z^2$, which we can pull over to an expectation $E$ on $C^*(H(\Z))$ such that
\[
E(\bar u_m\bar v^k)=
\begin{cases}
0&\text{unless $k=0$ and $m_2=0$}\\
\bar u_{m_1,0}&\text{if $k=0$ and $m_2=0$.}
\end{cases}
\]
A direct calculation shows that $E$ has the tracial property $E(ab)=E(ba)$ ($E$ is a \emph{centre-valued trace} on $C^*(H(\Z))$\;), and the isomorphism of $C^*(H(\Z))$ onto $C(\T)\rtimes_{\id,\omega}\Z^2$ carries $\sigma$ into $t\mapsto \widehat\beta_{(e^{it},1)}$. Thus any state of the form $\phi\circ E$ is an invariant trace of $C^*(H(\Z))$. These are the traces described in Example~\ref{tracesHgp}.
\end{remark}

\end{document}